\documentclass[12pt,twoside]{article} 
\usepackage{amsmath,amssymb,mathrsfs,graphicx}
\usepackage[small,sc]{caption} 

\usepackage[nodvipsnames]{color}

\newtheorem{theorem}{Theorem}[section]
\newtheorem{lemma}[theorem]{Lemma}
\newtheorem{proposition}[theorem]{Proposition}

\newtheorem{definition}[theorem]{Definition}

\newtheorem{rmrk}[theorem]{Remark}

\DeclareMathAlphabet{\mathbfit}{OML}{cmm}{b}{it}

\makeatletter
\@addtoreset{equation}{section}
\makeatother

\newcommand{\fig}[3] {
\begin{figure}[htb]
  \centering
  \vspace*{-10pt}         
  \includegraphics[width=#2]{imifp-#1.pdf}
  \begin{minipage}[t]{0.80\linewidth} 
    \vspace*{-24pt}        
    \caption{#3}
    \protect\label{#1}
  \end{minipage}
\end{figure}
\smallskip
}

\newenvironment{remark}
{\begin{rmrk} \em}
{\end{rmrk}}

\newcommand{\fn} {function}
\newcommand{\me} {measure}

\newcommand{\erg} {ergodic}
\newcommand{\sy} {system}

\newcommand{\pr} {probability}
\newcommand{\dsy} {dynamical system}
\renewcommand{\o} {orbit}

\newcommand{\R} {\mathbb{R}}
\newcommand{\C} {\mathbb{C}}
\newcommand{\Q} {\mathbb{Q}}
\newcommand{\Z} {\mathbb{Z}}
\newcommand{\N} {\mathbb{N}}
\newcommand{\qed} {\hfill {\small Q.E.D.} \par\medskip}
\newcommand{\skippar} {\par\medskip}
\newcommand{\ds} {\displaystyle}
\newcommand{\proof} {\noindent \textsc{Proof.} }
\newcommand{\proofof}[1] {\noindent \textsc{Proof of {#1}.} }
\newcommand{\article}[3] {\textsc{{#1}}, {\itshape {#2}}, {{#3}}.}
\newcommand{\book}[3] {\textsc{{#1}}, {\itshape {#2}}, {{#3}}.}
\newcommand{\vol} {\textbf}
\newcommand{\eps} {\varepsilon}

\newcommand{\rset}[2] {\left\{ #1 \: \left| \: #2 \right. \! \right\} }

\newcommand{\symmdiff} {\triangle}

\newcommand{\into} {\longrightarrow}

\renewcommand{\emptyset} {\varnothing}

\newcommand{\m} {mixing}
\newcommand{\ob} {observable}

\newcommand{\ui} {(0,1)}   
\newcommand{\rp} {\R^+}   
\newcommand{\br} {\tau}   
\newcommand{\ibr} {\phi} 
\newcommand{\pa} {\mathscr{P}}   
\newcommand{\ind} {\mathcal{J}}   
\newcommand{\leb} {m}   
\newcommand{\bj} {\mathbfit{j}}  
\newcommand{\sca} {\mathscr{A}}   
\newcommand{\avg} {\overline{\mu}}   
\newcommand{\avgleb} {\overline{\leb}}   
\newcommand{\scv} {\mathscr{V}}
\newcommand{\ivui} {V \nearrow \ui}
\newcommand{\ivrp} {V \nearrow \rp}

\newcommand{\go} {\mathcal{G}}   
\newcommand{\lo} {\mathcal{L}}   
\renewcommand{\P} {\mathbb{P}}
\newcommand{\pippo} {{p=0} \atop {p\equiv j\, (\text{mod } q)}}
\renewcommand{\aa} {\chi}
\newcommand{\bb} {\psi}
\renewcommand{\t} {\theta}   
\newcommand{\rat}{\mathcal{M}}   

\begin{document}

\title{\textbf{Infinite mixing for one-dimensional maps with an indifferent 
fixed point}}

\author{
\scshape
Claudio Bonanno\,\thanks{
Dipartimento di Matematica, Universit\`a di Pisa, Largo Bruno 
Pontecorvo 5, 56127 Pisa, Italy. E-mail: 
\texttt{claudio.bonanno@unipi.it}.}
,
Paolo Giulietti\,\thanks{Centro di Ricerca Matematica 
``Ennio de Giorgi'', Scuola Normale Superiore, Piazza dei 
Ca\-valieri 7, 56126 Pisa, Italy.
E-mail: \texttt{paolo.giulietti@sns.it}.}
,
Marco Lenci\,\thanks{
Dipartimento di Matematica, Universit\`a di Bologna,
Piazza di Porta San Donato 5, 40126 Bologna, Italy. 
E-mail: \texttt{marco.lenci@unibo.it}.}
\thanks{
Istituto Nazionale di Fisica Nucleare,
Sezione di Bologna, Via Irnerio 46,
40126 Bologna, Italy.}
}

\date{Final version for \emph{Nonlinearity} \\[6pt]
September 2018}

\maketitle

\begin{abstract}
  We study the properties of `infinite-volume mixing' for two classes of 
  intermittent maps: expanding maps $[0,1] \longrightarrow [0,1]$ with an 
  indifferent fixed point at 0 preserving an infinite, absolutely continuous 
  measure, and expanding maps $\R^+ \longrightarrow \R^+$ with an 
  indifferent fixed point at $+\infty$ preserving the Lebesgue measure. 
  All maps have full branches. While certain properties are 
  easily adjudicated, the so-called global-local mixing, namely the 
  decorrelation of a global and a local observable, is harder to prove. We 
  do this for two subclasses of systems. The first subclass includes, 
  among others, the Farey map. The second class includes the standard 
  Pomeau-Manneville map $x \mapsto x+x^2$ mod 1. Morevoer, we 
  use global-local mixing to prove certain limit theorems for our 
  intermittent maps.

\end{abstract}

\section{Introduction}
\label{sec-intro}

Expanding maps of the interval with indifferent, a.k.a.\
neutral, fixed points
are among the most intensively studied classes of \dsy s. 
They are considered the easiest examples of non-uniformly 
hyperbolic maps, where the mechanism that induces chaoticity 
is not as favorable---and somehow special---as in uniformly hyperbolic 
maps.

An indifferent fixed point can dramatically change the 
dynamical properties of an otherwise uniformly expanding map. 
Trajectories will spend long stretches of time in a neighborhood of the 
fixed point, nearly motionless, before returning to the strongly
expanding region of the space, where they exhibit a seemingly random
motion. In the physical literature, this behavior has been called 
\emph{intermittence}, and maps with indifferent fixed points sometimes
referred to as \emph{intermittent maps}. They have been widely used 
as models for a variety of ``anomalous'' dynamical phenomena. 
A representative, far from exhaustive, list of references includes 
\cite{gt, gnz, bg, zk, k}.

If the fixed point is \emph{strongly neutral}, which means that
the second derivative is continuous there, these \sy s preserve 
a Lebesgue-absolutely continuous infinite \me\ under very
general conditions \cite{t80}.
This and the fact that uniformly expanding interval maps are standard 
and somewhat elementary \dsy s has led to intermittent maps of the 
interval being very popular in the field of infinite \erg\ theory 
\cite{t80, t83, a, t00, z, i03}, considering also the many applications 
of its most notorious example, the Farey map \cite{d, p, i, ks, he, kms}.

Here we are interested in their \m\ properties, especially in the sense of
the recent definitions of \emph{infinite \m} given by Lenci 
\cite{limix, lpmu}. The expression `infinite \m' refers to all the 
notions, or formal definitions, which are supposed to replace, or extend, 
the definition of \m\ of finite \erg\ theory. 

The quest for an effective notion of infinite \m\ has a long history (a 
short version of which may be found in the introduction of \cite{limix}). 
Recent times have seen a significant surge of interest in this subject, 
both on its foundational aspects and on the application of new, 
sophisticated techniques to old problems \cite{limix, dr, liutam, mt12, 
ko, a13, lpmu, te, rt, lt, a16, mt17, lmmaps, s, dn}.

In \cite{mt12, te} Melbourne and Terhesiu studied a large class of 
interval maps with an indifferent fixed point, obtaining 
strong results related to the notion of \m\ first envisaged by Hopf 
in 1937 \cite{h} and later formalized, in slightly different ways, by 
Krickeberg \cite{kr}, Papangelou \cite{pa} and Friedman \cite{fr}.
This notion is sometimes referred to as \emph{Hopf-Krickeberg \m} or
\emph{rational \m}. In the case of a map $T$ preserving an infinite 
\me\ $\mu$, it  corresponds to the existence of a \emph{scaling rate} 
$(\rho_n)_{n \in \N}$ such that
\begin{equation}
  \label{def-rat-m}
  \lim_{n \to \infty} \rho_n \, \mu( (f \circ T^n) g ) = \mu(f) \mu(g),
\end{equation}
for all $f,g$ in certain subspaces of $L^1(\mu)$. Here, as usual, 
$\mu(f)$ is short for $\int f \, d\mu$. (See also \cite{a13} for 
the definition of \emph{rational weak \m}.)

From the point of view of the stochastic properties of \dsy s,
(\ref{def-rat-m}) corresponds to a local limit theorem. In the 
terminology used in the present paper, it represents a strong form 
of \emph{local-local \m}; cf.\ Section \ref{subs-im-ui}.

\skippar

The definitions of infinite \m\ introduced in \cite{limix}---also
referred to as \emph{infinite-volume \m}---hinge on the concept 
of a \emph{global \ob}. Informally speaking, a global 
\ob\ is a bounded \fn\ that is supported more or less throughout the
phase space, as opposed to a \emph{local \ob}, which is akin to a 
compactly supported \fn. In the present context, if $T: [0,1] \into 
[0,1]$ has a neutral fixed point at 0 and preserves an infinite \me\ 
$\mu$ which assigns finite mass to all $[a, 1]$, 
a global \ob\ is any $F \in L^\infty(\mu)$ for which
\begin{equation}
  \label{def-avg-intro}
  \avg(F) := \lim_{a \to 0^+} \frac1 {\mu([a,1))} \int_a^1 F \, d\mu 
\end{equation}
exists. In other words, a global \ob\ is a bounded \fn\ $F$ whose 
averages over larger and larger portions of the space (in the sense of 
the \me) converge to an \emph{infinite-volume average} $\avg(F)$.
A local \ob\ is any $g \in L^1(\mu)$. (For the sake of readability,
global and local \ob s are indicated, respectively, with uppercase and
lowercase letters.)

We speak of \emph{global-global \m} when, for every pair of 
global \ob s $F,G$,
\begin{equation}
  \label{ggm1-intro}
  \lim_{n \to \infty} \avg( (F \circ T^n) G ) = \avg(F) \avg(G).
\end{equation}
We call \emph{global-local \m} the case when, for all global \ob s $F$ 
and local \ob s $g$,
\begin{equation}
  \label{glm2-intro}
  \lim_{n \to \infty} \mu( (F \circ T^n) g ) = \avg(F) \mu(g).
\end{equation}
Both definitions have other versions as well, which are discussed in 
Section \ref{subs-im-ui}. 

In \cite{limix} global-global \m\ and global-local \m\ were proved, 
under suitable conditions, for \dsy s representing random walks in 
$\Z^d$. In \cite{lmmaps} both types of \m\ were verified for a certain 
class of uniformly expanding maps of the real line, the so-called 
\emph{quasi-lifts} and their local modifications.

In this paper we show how maps with an indifferent fixed point
of the type outlined earlier can never be global-global \m, and
present a general method to prove global-local \m\ for such \sy s. 
The method covers a large class of examples, including the Farey 
map and many Pomeau-Manneville maps. It becomes particularly 
simple when, via conjugation, we represent our maps as \dsy s on 
$\R^+$ preserving the Lebesgue \me.

To summarize, the various sections of the paper are organized as 
follows. Section \ref{sec-setup} is the backbone: we describe  
in detail the classes of maps we study, review the notions of global 
and local \ob s, together with the various definitions of infinite \m, and
state our results. In Section \ref{sec-appl} we present two 
examples of limit theorems that can be proved for intermittent maps 
that are global-local \m. The rest of the paper is devoted to the proofs.
In Section \ref{sec-pfs-easy} we prove the simpler results.
In Section \ref{sec-glm2} we give the scheme of the proof of 
global-local \m. This is based on the existence of a local \ob\ with 
certain monotonicity properties. Such existence will be established, 
for all cases considered,  in Section \ref{sec-pers-mon}. The proof 
of global-local \m\ also uses the exactness of the map, which is a 
standard result. However, for the class of maps $\rp \into \rp$ that 
we study, we found no proof in the literature, so we give our own 
proof in Appendix \ref{sec-ex}. Finally, Appendix \ref{sec-tech} 
contains the proofs of two technical results.

\bigskip\noindent
\textbf{Acknowledgments.}\ This research is part of the authors' 
activity within the \linebreak DinAmicI community, see 
\texttt{www.dinamici.org}. C.~Bonanno and M.~Lenci thank 
the Istituto Nazionale di Alta Matematica and its division 
Gruppo Nazionale di Fisica Matematica for various forms of 
support. P.~Giulietti thanks the Universidade Federal do Rio 
Grande do Sul, Porto Alegre, Brazil, where part of this work was
done. He also acknowledges the financial support of 
the Centro di Ricerca Matematica ``Ennio de Giorgi'' and of 
UniCredit Bank R\&D Group, through the ``Dynamics and 
Information Theory Institute'' at the Scuola Normale Superiore.

\section{Setup and results}
\label{sec-setup}

In this section we give a detailed presentation of the maps we 
consider in the paper. We divide them in two classes: maps 
$[0,1] \into [0,1]$ with a strongly neutral fixed point at 0, and maps 
$\rp \into \rp$ with a neutral fixed point at $+\infty$ preserving the
Lebesgue \me, cf.\ Figs.\
\ref{alpha03} and \ref{rp3br} later in the section. They are 
morally the same \sy s, because one can always pass from one type 
of map to the other via a suitable conjugation. But the conjugation 
will not map the first class exactly onto the second, hence the need to 
distinguish the two cases.

In order to emphasize the similar nature of the maps in the two
classes, we choose to always use an open phase space. This 
means that for the rest the paper `unit interval' will always indicate 
the open interval $\ui$. This choice has no consequence on our 
results. 

\subsection{Maps of the unit interval}
\label{subs-ui}

In the case of a map $T : \ui \into \ui$, we
assume there to be a finite or infinite sequence
of numbers $0 = a_0 < a_1 < \ldots < a_k < \ldots \le 1$.
If the sequence is finite, its last element is $a_N := 1$; 
in this case we set $\ind := \{ 0, \ldots, N-1 \}$. If the sequence is 
infinite, $\lim_n a_n = 1$; in this case we
set $\ind = \N$ (in our notation $0 \in \N$). For $j \in \ind$, 
denote $I_j := (a_j, a_{j+1})$. Thus, $\pa := \{ I_j \}_{j \in \ind}$
is a partition of $\ui$ mod $\leb$, the Lebesgue \me\ on $\R$.

We assume that $T$ is a Markov map w.r.t.\ $\pa$, with the
following properties:
\begin{itemize} 
\item[(A1)] $T |_{I_j}$ possesses an extension $\br_j : [a_j, a_{j+1}]
  \into [0,1]$ which is bijective and $C^2$ up to the boundary.
  
\item[(A2)] There exists $\Lambda > 1$ such that $| \br'_j | \ge
  \Lambda$, for all $j \ge 1$.
  
\item[(A3)] There exists $K > 0$ such that $\ds \frac{ |\br''_j | }
  { |\br'_j |^2 } \le K$, for all $j\ge 0$.
  
\item[(A4)] $\br_0$ is convex with $\br_0(0) = 0$, $\br'_0(0) = 1$,
  and $\br'_0(x) > 1$, for $x \in (0,a_1]$.
\end{itemize}

The following statements, which were proved, respectively, in 
\cite{t80} and \cite{t83}, will be useful in the remainder. 

\begin{theorem} \label{thm-thaler} 
  Under the assumptions {\rm (A1)-(A4)},
  \begin{itemize}
  \item[(a)] $T$ preserves an infinite invariant \me\ $\mu$ which is 
    absolutely continuous w.r.t.\ the Lebesgue \me\ $\leb$ and, up
    to multiplicative constants, is the unique absolutely continuous
    invariant \me. Moreover, the infinite density $h := d\mu / d\leb$ 
    is positive and unbounded only near 0.
    
  \item[(b)] $T$ is conservative and exact (w.r.t.\ $\leb$ or $\mu$, which 
    is the same).
  \end{itemize}
\end{theorem}

We recall that $T$ is said to be exact when, denoted by $\sca$ the 
$\sigma$-algebra of its reference space, the \emph{tail $\sigma$-algebra} 
$\bigcap_{n \in \N} T^{-n} \sca$ is trivial, i.e., it contains only null sets or 
complements of null sets.

Exactness is a strong \m\ property which has the distinct advantage
of being defined in the same way both in finite and infinite
\erg\ theory. Within the scope of the present paper, it has the
additional merit of being a key ingredient for the proof of 
the global-local \m\ (\ref{glm2-intro}).

\subsection{Infinite mixing for maps of the unit interval}
\label{subs-im-ui}

For $F \in L^\infty(\ui, \mu)$ and $a \in (0,1)$, denote
\begin{align}
\label{mu-v-ui}
  \mu_{[a,1)}(F) &:= \frac1 {\mu([a,1))} \int_a^1 F \, d\mu ; \\[4pt]
\label{def-avg-ui}
  \avg(F) &:= \lim_{a \to 0^+} \mu_{[a,1)}(F).
\end{align}
The limit (\ref{def-avg-ui}) might not exist. When it does,
we say that $F$ is a \emph{global \ob} and call $\avg(F)$ the
\emph{infinite-volume average} of $F$. The space of all global 
\ob s is denoted by $\go$. In addition, we call any $f \in \lo := 
L^1(\ui, \mu)$ a \emph{local \ob}. 

\begin{remark}  
  In the framework of \cite{limix} and \cite{lpmu} the definitions
  (\ref{mu-v-ui})-(\ref{def-avg-ui}) correspond to choosing the
  \emph{exhaustive family} $\scv := \rset{[a,1)} {0 < a < 1}$. An
  exhaustive family is a collection of finite-\me\ sets that play the
  role of ``large boxes'' in a reference space. The generic element  
  of $\scv$ will also be denoted $V$. The limit $a \to 0^+$ is called 
  the \emph{infinite-volume limit}. In more suggestive notation we will
  also indicate it by $\ivui$.
\end{remark}

To visualize an example of a global \ob, one can think of a bounded
\fn\ of $(0,1)$ which oscillates around 0 in such a way that the
limit in (\ref{def-avg-ui}) exists. A more intuitive visualization
of a global \ob s will be given in the Section \ref{subs-rp}, where
the reference space is $\rp$. Notice that a bounded \fn\
which has a limit at 0 is also a global \ob, but a very
insignificant one, because it is arbitrarily close to a constant in all but 
a tiny fraction of the space (in the sense of the \me). Unquestionably,
any reasonable definition of mixing must be trivially verified on constant 
\ob s.

We briefly recall the definitions of `infinite-volume \m' presented in 
\cite{limix, lpmu}. The \dsy\ $(\ui, \mu, T)$ is called 
\emph{global-local \m} of type
\begin{description}
\item[(GLM1)] \, if, $\forall F \in \go$, $\forall g \in \lo$ with
  $\ds \mu(g)=0$, $\ds \lim_{n \to \infty} \, \mu((F \circ T^n) g) =
  0$;
  
\item[(GLM2)] \, if, $\forall F \in \go$, $\forall g \in \lo$, $\ds
  \lim_{n \to \infty} \, \mu((F \circ T^n) g) = \avg(F) \mu(g)$;

\item[(GLM3)] \, if, $\forall F \in \go$, $\ds \lim_{n \to \infty} \,
  \sup_{g \in \lo \setminus 0 } \frac{ \left| \mu((F \circ
  T^n) g) - \avg(F) \mu(g) \right| } {\mu( |g| )}= 0$.
\end{description}
It is called called \emph{global-global \m} of type
\begin{description}
\item[(GGM1)] \, if, $\forall F,G \in \go$, $\ds \lim_{n \to \infty}
  \, \avg ( (F \circ T^n) G ) = \avg(F) \, \avg(G)$;

\item[(GGM2)] \, if, $\forall F,G \in \go$, $\ds \lim_{{\ivui} \atop {n \to 
  \infty}} \mu_V ((F \circ T^n) G) = \avg(F) \, \avg(G)$.
\end{description}
The limit in {\bf (GGM2)} means that, for all $\eps>0$, there exists 
$M>0$ such that the l.h.s., defined as in (\ref{mu-v-ui}), is
$\eps$-close to the limit for all $V = [a,1)$, with $\mu(V) \ge M$, and
all $n \ge M$. It is called the `joint infinite-volume and time limit'; cf.\ 
\cite[Defn.~2.2]{lpmu}. 

Our first proposition states that if $T$ is such that (A1)-(A4) are 
satisfied, then $\avg$ is an invariant functional for the dynamics. If
this were not the case, the above definitions would not make sense. 
To keep the exposition fluid, we postpone the proof to Section 
\ref{sec-pfs-easy}.

\begin{proposition} \label{prop-inv-ui}
  Let $T: \ui \into \ui$ verify {\rm (A1)-(A4)}. For all $F \in \go$ and $n \in \N$,
  $\avg(F \circ T^n)$ exists and equals $\avg(F)$.
\end{proposition}

Finally, the \sy\ is called \emph{local-local \m}
\begin{description}
\item[(LLM)] \, if, $\forall f \in \lo \cap \go$, $g \in \lo$,
  $\ds \lim_{n \to \infty} \, \mu ( (f \circ T^n) g ) = 0$.
\end{description}
Since, in the present case, $\go$ comprises all $F \in L^\infty$ which 
possess an infinite-volume average and  $\lo = L^1$, one verifies that 
{\bf (LLM)} is equivalent to the definition of zero-type \dsy:
$\lim_{n \to \infty} \mu(T^{-n} A \cap A) = 0$, for all $A$ with $\mu(A)
< \infty$ \cite{hk, ds}.

The property with which we are most concerned in this article 
is {\bf (GLM2)}, which can be recast like this: For every 
$\mu$-absolutely continuous \pr\ \me\ $\nu$ and every global \ob\ $F$,
\begin{equation}
  \lim_{n \to \infty} T^n_* \nu (F) = \avg(F),
\end{equation}
where $T^n_* \nu = \nu \circ T^{-n}$ denotes the push-forward
of $\nu$ via the map $T^n$. In this sense, {\bf (GLM2)}
represents a very weak form of convergence of $T^n_* \nu$ to
$\mu$, which cannot occur in any conventional sense, as the former
are \pr\ \me s and the latter is an infinite \me.

For all the other properties we have the following.

\begin{proposition} \label{prop-other-ui}
  A map $T: \ui \into \ui$ verifying {\rm (A1)-(A4)} is {\bf (GLM1)} and 
  {\bf (LLM)}, but not {\bf (GLM3)}, {\bf (GGM1)} or {\bf (GGM2)}.
\end{proposition}

Once again, we give the proof of Proposition \ref{prop-other-ui} in 
Section \ref{sec-pfs-easy}. We will see that, in the present case, it is 
relatively easy to check all of the conditions except {\bf (GLM2)}. This 
does not mean, however, that these definitions are unimportant or give 
no information about the \sy; quite the contrary. For example, the fact
that $T$ cannot be global-global \m\ formalizes the idea that an 
expanding map with an indifferent fixed point has radically different
chaotic properties than a uniformly expanding map. This is no surprise, 
given that the former is very close to the identity in the overwhelming 
majority of the space (in terms of the \me). By way of comparison, we 
observe that the uniformly expanding maps on $\R$ studied in 
\cite{lmmaps} are generally expected to be global-global \m.
(For further comparison with the results of \cite{lmmaps} see
the last paragraph of Section \ref{sec-glm2}.)

We now introduce a class of maps satisfying (A1)-(A4) which verify 
{\bf (GLM2)}. They are Markov maps with $N=2$ surjective branches.
(The case $2 < N < \infty$ can be treated as well, though the 
necessary hypotheses become more cumbersome, cf.\ Remark
\ref{rk-claudio} below.) 

In view of (A1), let us denote $\ibr_0 := \br_0^{-1} : [0,1] \into [0,a_1]$ 
and $\ibr_1 := \br_1^{-1} : [0,1] \into [a_1,1]$. These \fn s, which extend
the inverse branches of $T$, are bijective and $C^2$ up to the boundary.
Moreover,  $\ibr_0'(0)=1$, $\ibr_0'(x) \in (0,1)$ for $x\in (0,1]$, and $\ibr_0$ 
is concave. Recalling that $h$ is the density of the infinite invariant 
measure $\mu$ given by Theorem \ref{thm-thaler}\emph{(a)}, we make 
the following extra assumptions:

\begin{itemize}
\item[(A5)] $\ibr_1$ is decreasing (equivalently, $\br_1$ is decreasing).

\item[(A6)] $\ibr_0 + \ibr_1$ is increasing and concave.

\item[(A7)] $\ibr'_0 \, (h\circ \ibr_0)/h$ is differentiable, strictly decreasing 
  and convex.

\item[(A8)] $\ibr_0'\, (h \circ \ibr_0) + \ibr_1'\, (h \circ \ibr_1) \ge 0$.
\end{itemize}

\begin{remark} \label{rk-a8}
  If $h$ is decreasing, (A8) follows from (A6). In fact, $h>0$, $\ibr'_0>0$ and 
  $\ibr_0 \le \ibr_1$ imply $\ibr_0'\, (h \circ \ibr_0) \ge \ibr_0'\, (h \circ \ibr_1) 
  \ge -\ibr_1'\, (h\circ \ibr_1)$. 
\end{remark}

\begin{theorem} \label{main-thm-ui} 
  Let $T: \ui \into \ui$ satisfy assumptions {\rm (A1)-(A8)} w.r.t.\ $\pa = 
  \{I_0, I_1\}$. Then $T$ is {\bf (GLM2)}.
\end{theorem}  

The proof of this theorem is given in Sections \ref{sec-glm2}
and \ref{sec-pers-mon}. An interesting family of maps which 
satisfy the hypotheses of the theorem is constructed starting from the 
Farey map:
\begin{equation} \label{farey}
  T_0(x) = \left\{ \begin{array}{ll}
  \ds \frac{x}{1-x}\, ,& x\in \left[ 0 , \frac12 \right]; \\[14pt]
  \ds \frac{1-x}{x}\, , & x\in \left[ \frac12, 1 \right].
  \end{array} \right.
\end{equation}
It is well known that $T_0$ preserves the infinite \me\ on $\ui$ whose
density is $h(x) = 1/x$. The inverse branches of $T_0$ are easily computed
to be $\ibr_0(x) = x/(1+x)$ and $\ibr_1(x) = 1/(1+x)$.

For $\alpha \in (0,1)$, set $I_0 := (0,2^{\alpha-1})$ and 
$I_1 := (2^{\alpha-1},1)$, and consider the map $T_\alpha: \ui \into \ui$ 
implicitly defined by the inverse branches
\begin{equation}
  \ibr_0(x) := \frac{x}{(1+x)^{1-\alpha}}, \qquad 
  \ibr_1(x) := \frac{1}{(1+x)^{1-\alpha}},
\end{equation}
where $\ibr_j : I_j \into [0,1]$, for $j \in \{0,1\}$. An example is shown in 
Fig.~\ref{alpha03}. We have
\begin{equation}
  \ibr'_0(x) = \frac{1+\alpha x}{(1+x)^{2-\alpha}}, \qquad 
  \ibr'_1(x) = - \frac{1-\alpha}{(1+x)^{2-\alpha}}
\end{equation}
and
\begin{equation}
  \ibr''_0(x) = - (1-\alpha)\, \frac{2+\alpha x}{(1+x)^{3-\alpha}}, \qquad 
  \ibr''_1(x) = \frac{(2-\alpha)(1-\alpha)}{(1+x)^{3-\alpha}}.
\end{equation}

\fig{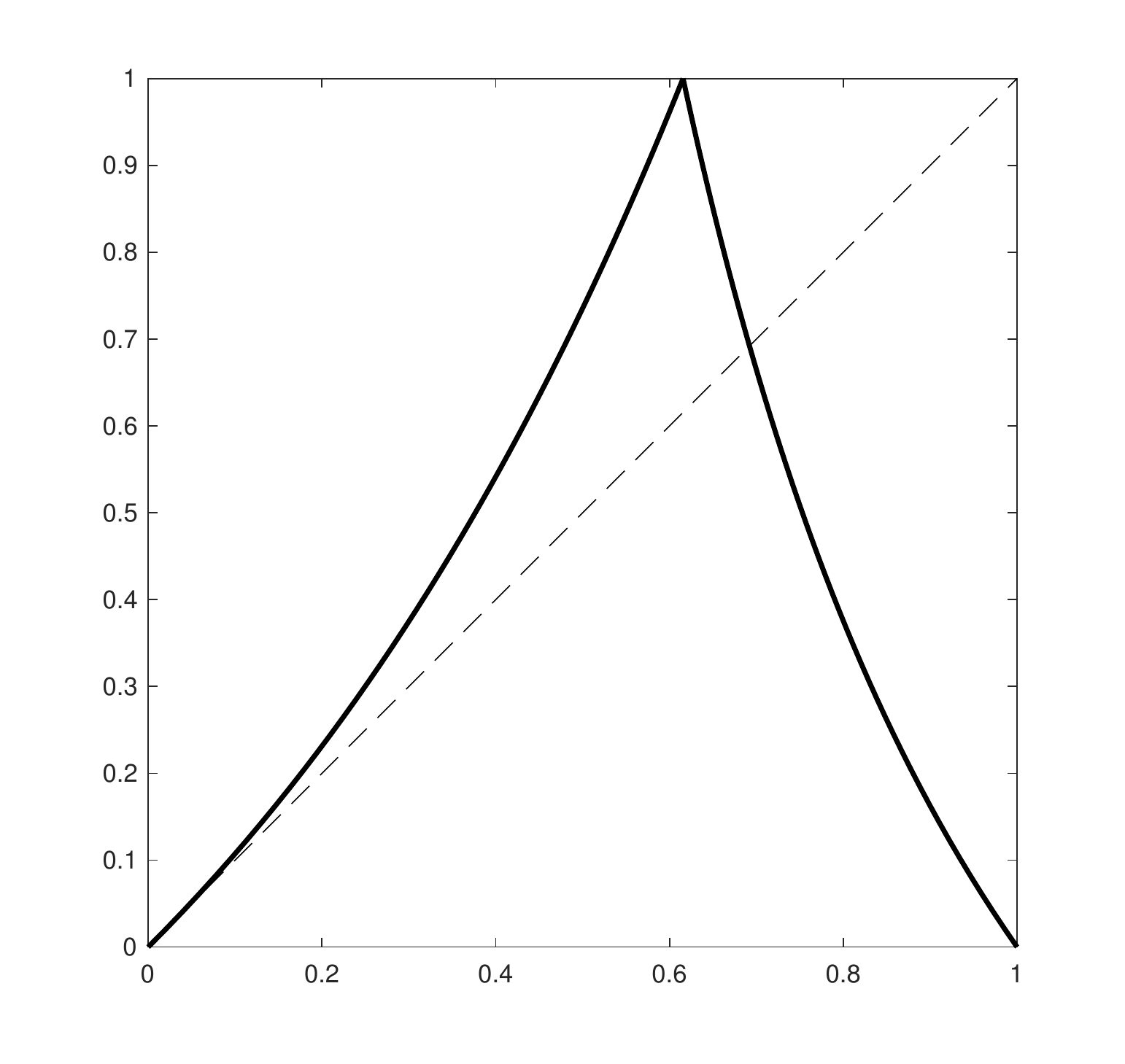}{9cm}{The map $T_\alpha$ defined in Section \ref{subs-im-ui},
for $\alpha = 0.3$.}

It is easy to check that $T_\alpha$ verifies (A1)-(A6). Moreover,
$T_\alpha$ preserves the same \me\ preserved by the 
Farey map $T_0$. In fact, given $h(x)=1/x$, one has
\begin{equation} \label{f-like40}
  |\ibr'_0|\, (h \circ \ibr_0) + |\ibr'_1|\, (h \circ \ibr_1) = h,
\end{equation}  
which implies that $\int (F \circ T_\alpha) h\, d\leb = \int F\, h\, d\leb$,
for all $F \in L^\infty$.
(In other words, if $P$ denotes the transfer operator of $T_\alpha$
relative to $\mu$, cf.\ (\ref{trans-op-ui}), the identity (\ref{f-like40})
is equivalent to $P1 = 1$, where $1$ is the (non-integrable) \fn\
which is identically equal to 1.) Finally, the equation
\begin{equation} 
  \ibr'_0(x) \, \frac{h(\ibr_0(x))}{h(x)} = \frac{1+\alpha x}{(1+x)} ,
\end{equation}  
proves (A7), while (A8) follows from (A6) and the monotonicity of 
$h$, as pointed out in Remark \ref{rk-a8}.

\begin{remark} \label{rk-farey}
  Since the parameter $\alpha$ ranges in $(0,1)$, the above family
  does \emph{not} include the Farey map (\ref{farey}). The problem is
  that $T_0'(1) = -1$ and (A2) is not verified. But the conclusions of 
  Theorem \ref{main-thm-ui} hold for the Farey map too. As it will
  be clear later, cf.\  Definition \ref{def-per-mon} and Theorem 
  \ref{thm-main}, it is sufficient to find a persistently monotonic local 
  \ob\ for $T_0$. This was done in \cite[Lem~8.13]{i}.
\end{remark}

\begin{remark} \label{rk-claudio}
  Theorem \ref{main-thm-ui} can be extended to the case of $N$
  branches, $2 < N < \infty$, if, in addition to (A1)-(A4), the following 
  assumptions are made:
  \begin{itemize}
    \item $\br_k$ is increasing and convex for all $k \in \{0,\dots,N-2\}$;
      $\br_{N-1}$ is decreasing.
    \item $\sum_{k=0}^{N-1} \ibr_k$ is increasing and concave.
    \item $\ibr'_k \, (h\circ \ibr_k)/h$ is strictly decreasing and convex 
      for all $k \in \{0,\dots,N-2\}$.
    \item $h$ is decreasing (or the analogue of (A8) holds with $\ibr_k$ 
      in place of $\ibr_0$, for all $k\in \{0,\dots,N-2\}$).
  \end{itemize}
  The proof of this generalization adds computations but no new ideas to
  the one presented in the paper, so we omit it. 
\end{remark}

A very popular family of intermittent maps of the unit interval is the 
loosely defined class that goes by the names of Pomeau and 
Manneville. These maps have been introduced to study in full rigor 
certain intermittency phenomena initially described by Pomeau and 
Manneville in the 1980's \cite{pm,m}. Although no precise definition 
exists, most mathematicians would agree that a Pomeau-Manneville 
map is a map of $\ui$ with two increasing branches satisfying at 
least (A1) and (A4). It is natural to ask weather maps of this type 
are global-local \m. Neither Theorem \ref{main-thm-ui} nor 
Remark \ref{rk-claudio} address this case because they assume
one branch to be decreasing. Nonetheless, many Pomeau-Manneville
maps are {\bf (GLM2)}. For a given \sy, this can be shown via the 
results of Section \ref{subs-im-rp} below, provided one has enough 
information about the invariant \me\ $\mu$. The problem, however,
is that the general theorems that are available at this time 
\cite{t80, t01} do not provide enough control on $\mu$. We refer the 
reader to Remarks \ref{rk-pm2} and \ref{rk-pm3} in Section \ref{subs-im-rp}.

\subsection{Maps of the half-line}
\label{subs-rp}

Given a map $T : \rp \into \rp$, we assume 
that there exists a finite or 
infinite sequence $a_1 > a_2 > \ldots > a_k > \ldots \ge 0$.
If the sequence is finite, its last element is $a_N := 0$; 
in this case $\ind := \{ 0, \ldots, N-1 \}$. If the sequence is 
infinite, $\lim_n a_n = 0$; in this case $\ind = \N$. Denote
$I_0 := (a_1, +\infty)$ and, for $j \in \ind \setminus \{0\}$,
$I_j := (a_{j+1}, a_j)$. Once again, $\pa := \{ I_j \}_{j \in \ind}$
is a partition of $\R^+$ mod $\leb$. 

We also assume that:

\begin{itemize} 
\item[(B1)] $T |_{I_j}$ is a bijective map onto $\R^+$,
  and possesses an extension $\br_j$ which, for $j=0$, is
  defined on $[a_1, +\infty)$ and, for $j \ge 1$, is defined on 
  $[a_{j+1}, a_j)$ or $(a_{j+1}, a_j]$. $\br_j$ is $C^2$ 
  up the boundary.
  
\item[(B2)] There exists $\Lambda > 1$ such that $| \br'_j | \ge
  \Lambda$, for all $j \ge 1$.
  
\item[(B3)] There exists $K > 0$ such that $\ds \frac{ |\br''_j | }
  { |\br'_j |^2 } \le K$, for all $j\ge 0$.
  
\item[(B4)] The \fn\ $u(x) := x - \br_0 (x)$ is positive, convex and 
  vanishing (hence decreasing), as $x \to +\infty$. Furthermore, 
  $u''$ is decreasing (hence vanishing).
  
\item[(B5)] $T$ preserves the Lebesgue \me\ $\leb$.
\end{itemize}

The most restrictive assumption here, compared to Section
\ref{subs-ui}, is (B5): we require $T$ to preserve not just
an absolutely continuous \me, but exactly the 
Lebesgue \me. Some of the results we obtain (for example, Theorem
\ref{thm-ex-leb} and Proposition \ref{prop-other-rp}) would also
hold in the case where $T$ preserves an absolutely continuous, 
infinite, locally finite \me. With assumption (B5), however, the
infinite-volume average of a global \ob\ is defined in a very natural
way, see (\ref{def-avgleb}).

Note that, given a $T_o: \ui \into \ui$ satisfying (A1)-(A4), it is
straightforward to find a conjugation $\Phi : \ui \into \rp$ such
that $T := \Phi \circ T_o \circ \Phi^{-1}$ verifies (B5). It suffices to
take $\Phi(x) := \int_x^1 h \, d\leb$, where $h$ is the Radon-Nikodym 
derivative mentioned in Theorem \ref{thm-thaler}\emph{(a)}. But $T$ 
might not verify the other assumptions. For instance, it might not be 
expanding.

In analogy with Theorem \ref{thm-thaler}\emph{(b)}, we have:

\begin{theorem} \label{thm-ex-leb} 
  Under assumptions {\rm (B1)-(B5)}, $T$ is conservative 
  and exact.
\end{theorem}

The proof of Theorem \ref{thm-ex-leb}---in fact, a generalization 
thereof---is given in Appendix \ref{sec-ex}.

The \ob s that we associate with these types of maps are completely
analogous to those defined in Section \ref{subs-im-ui}, with the difference
that use $\leb$ instead of $\mu$. More precisely, the class of
global \ob s is the space $\go$ of all $F \in L^\infty(\rp, \leb)$ such that
\begin{equation} \label{def-avgleb}
  \exists \avgleb(F) := \lim_{a \to +\infty} \leb_{(0,a]}(F) := 
  \lim_{a \to +\infty} \frac1a \int_0^a F \, d\leb.
\end{equation}
Correspondingly, 
the generic large box in reference space is $V = (0,a]$, and
the infinite-volume limit, here denoted $\ivrp$, is the limit
$a \to +\infty$. Finally, the class of local \ob s is $\lo := L^1(\rp, \leb)$. 

It is easy to see that any bounded periodic $F$ is a global \ob,
and $\avgleb(F)$ is the average of $F$ over a period.
Also, a large variety of ``quasi-periodic'' \fn s belong in $\go$, for 
instance $F(x) := e^{2\pi i x/\alpha} G(x)$, where $G$ is 
a bounded periodic \fn\ (in this case, if the ratio between $\alpha$ 
and the period of $G$ is irrational, $\avgleb(F) = 0$; otherwise
$F$ is periodic). More ``random'' \fn s also belong in $\go$: for 
example, if $f: \R \into \C$ is bounded and supported in 
$(0,b)$, and $(c_k)_{k \in \N}$ is a bounded sequence which 
possesses a Cesaro average, then $F(x) := \sum_{k\in \N} c_k 
f(x - kb)$ is a global \ob. 

\subsection{Infinite mixing for maps of the half line}
\label{subs-im-rp}

For $T: \rp \into \rp$ we consider the same definitions of 
infinite-volume \m\ presented in Section \ref{subs-im-ui}, with the
understanding that $\go$ and $\lo$ are those defined earlier,
$\mu$ is the Lebegue \me\ $\leb$, the exhaustive family is $\scv := 
\rset{(0,a]} {a > 0}$, and the infinite-volume limit is $\ivrp$ or, in
other words, $a \to +\infty$. The same results as in Section 
\ref{subs-im-ui} hold here, and they are again proved in Section
\ref{sec-pfs-easy}.

\begin{proposition} \label{prop-inv-rp}
  Let $T: \rp \into \rp$ verify {\rm (B1)-(B5)}. For all $F \in \go$ and 
  $n \in \N$, $\avgleb(F \circ T^n)$ exists and equals $\avgleb(F)$.
\end{proposition}

\begin{proposition} \label{prop-other-rp}
  A map $T: \rp \into \rp$ verifying {\rm (B1)-(B5)} is {\bf (GLM1)} and 
  {\bf (LLM)}, but not {\bf (GLM3)}, {\bf (GGM1)} or {\bf (GGM2)}.
\end{proposition}

We now introduce a class of maps satisfying (B1)-(B5) which verify 
{\bf (GLM2)}. They will be determined by the extra assumption:

\begin{itemize}
\item[(B6)] $\br_j$ is increasing and convex for all $j \ge 1$.
\end{itemize}

\fig{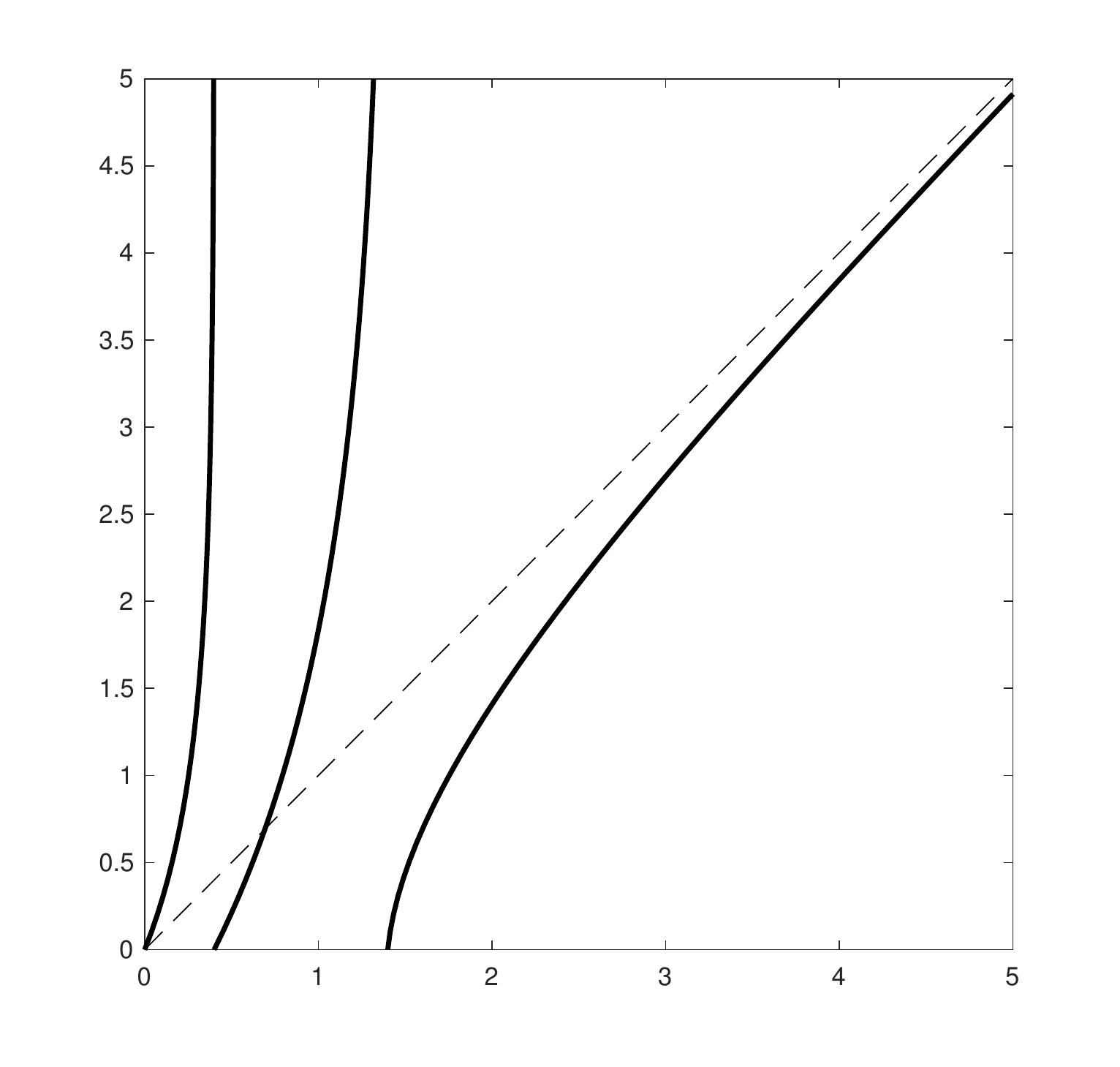}{9cm}{An example of a map $\rp \into \rp$ verifying (B1)-(B6).}

An example of such a map is shown in Fig.~\ref{rp3br}.
Once again, let $\ibr_j$ denote the inverse of $\br_j$. By (B1) and (B4), the 
functions $\ibr_j$ are bijective and $C^2$ up to the boundary, and 
$\ibr_0$ is increasing and convex. By (B6), $\ibr_j$ is increasing and
concave for all $j\ge 1$.

\begin{remark} \label{rk-conv-b6}
  If $T$ has only two branches, the convexity of $\br_1$ is a consequence
  of the other hypotheses. In fact, cf.\ the proof of Theorem \ref{exis-pm-rp}, 
  the preservation of the Lebesgue \me\ reads $\ibr_0' + \ibr_1' = 1$, whence
  $\ibr_0'' + \ibr_1'' = 0$. Therefore $\ibr_0$ and $\ibr_1$ have opposite 
  convexities. The same then holds for $\br_0$ and $\br_1$.
\end{remark}

\begin{theorem} \label{main-thm-rp} 
  Let $T: \rp \into \rp$ satisfy assumptions {\rm (B1)-(B6)}. Then $T$ 
  is {\bf (GLM2)}.
\end{theorem}

\begin{remark}
  The above theorem can be improved to include maps $T: \rp \into \rp$ 
  which verify (B1)-(B5) and have two branches, with $\br_1$ decreasing.
  In this case, however, we need to assume that the functions $\br_j$ are 
  $C^3$ up to the boundary, and add the following hypotheses:
  \begin{itemize}
  \item $\ibr_0 + \ibr_1$ is increasing.
  \item $\ibr''_1 - (\ibr'_1)^2 >0$.
  \item For each $x\in \rp$, one of the following two conditions holds: either
  \begin{displaymath}
    \ibr'''_1(x) + \ibr''_1(x) > 0 \quad \text{and} \quad 
    3 \ibr''_1(x) - (\ibr'_1(x))^2 + \ibr'_1(x) > 0 ;
  \end{displaymath}
  or
  \begin{displaymath}
    \ibr'''_1(x) + \ibr''_1(x) > (\ibr'_1(x))^2 .
  \end{displaymath}
  \end{itemize}
  We omit the proof of this extension for the same reasons as in Remark 
  \ref{rk-claudio}. The interested reader can nevertheless find it in 
  \cite{bgl}, where it is used to prove that the Boole map is {\bf (GLM2)}.
\end{remark}

\begin{remark} \label{rk-pm2}
  As discussed in Section \ref{subs-rp}, any map $T_o: \ui \into \ui$ 
  satisfying (A1)-(A4) can always be conjugated to a map 
  $T: \rp \into \rp$ that satisfies (B5), that is, preserves the Lebesgue \me. 
  The conjugation is $T := \Phi \circ T_o \circ \Phi^{-1}$, where 
  $\Phi(x) := \int_x^1 h \, d\leb$ and $h$ is an (infinite) invariant density 
  for $T_o$. Therefore, with enough information about $h$, one will be 
  able to check whether $T$ satisfies the hypotheses of Theorem 
  \ref{main-thm-rp}, thus showing that $T$ and therefore $T_o$ are 
  {\bf (GLM2)}. (Notice that all definitions of infinite-volume \m\ are 
  invariant w.r.t.\ the above conjugation.) For example, the 
  Pomeau-Manneville map $T_o(x) := x+x^2$ mod 1 has an 
  invariant density $h(x) = 1/x + 1/(1+x)$ \cite{t01}. By
  construction, the corresponding $T$ has two increasing full branches, 
  which are $C^\infty$ because $h$ is. It is a simple calculation to show 
  that the branch $\br_0$ is concave, implying (B6) via Remark
  \ref{rk-conv-b6}. The hypotheses (B2)-(B4) are also satisfied. Actually,
  as it will be clear later on (cf.\ Theorems \ref{thm-main} and 
  \ref{exis-pm-rp}), (B2)-(B4) are not directly needed in the proof of 
  global-local \m: they are only used to show that $T$ is exact, which 
  is another ingredient of the proof. But if one knows a priori that $T$ is 
  exact, which in this case follows from the exactness of $T_o$ 
  (Theorem \ref{thm-thaler}\emph{(b)}), one only need check 
  condition (B6) for $T$. For a map $T_o$ with two increasing 
  branches, by Remark \ref{rk-conv-b6}, this amounts to checking that
  \begin{equation}
    T'(x) = T_o'(\Phi^{-1}(x)) \,\frac{ h( T_o( \Phi^{-1}(x))) } 
    { h (\Phi^{-1}(x)) }
  \end{equation}
  is decreasing on the support of the branch $\br_0$, i.e, for $x > a_1$. 
  This is equivalent to 
  \begin{equation}
    \log T'(\Phi(x)) = \log T_o'(x) + \log h (T_o(x)) - \log h (x)
  \end{equation}
  being increasing for $0 < x < a_{o1}$, where $a_{o1} := \Phi^{-1}(a_1)$ 
  is the point of $\ui$ that separates the two Markov intervals of $T_o$.
\end{remark}

\begin{remark} \label{rk-pm3}
  The procedure that we have outlined in the previous remark can
  also be used in the opposite direction: starting from a map
  $T: \rp \into \rp$ verifying (B1)-(B6), and therefore {\bf (GLM2)}, one
  can construct an endless number of global-local \m\ maps
  $T_o := \Phi^{-1} \circ T \circ \Phi : \ui \into \ui$. This is in fact very 
  easy, as assumptions (B1)-(B6) are rather general and one 
  has ample freedom to choose $h$, and thus $\Phi$. For example,
  given a two-branched map $T$ satisfying (B1)-(B6), and chosen
  $h(x) := 1/x^p$, with $p \ge 1$, simple computations show that
  $T_o$ has a strongly neutral fixed point at 0 and is expanding
  and convex in a neighborhood of 0. That $T_o$ has two
  full increasing branches and preserves $\mu$, with
  $h = d\mu / d\leb$, is obvious by construction. So $T_o$ is  
  at least similar to a Pomeau-Manneville map. (In practice, most
  examples one will cook up lead to a fully convex branch at the origin.) 
  This observation shows that one can construct many global-local 
  \m\ maps of the unit interval with any index $p \ge 1$;  by this we 
  mean here that the invariant density has a singularity of the type 
  $x^{-p}$, as $x \to 0^+$.
  So it makes sense to study problems of linear response in this 
  context as well---e.g., the response of $\avg_p$ w.r.t.\ $p$
  (see \cite{bt, gg, brs} and references therein for the corresponding
  problem in the finite-\me\ case, possibly with noise).
\end{remark}

\section{Applications}
\label{sec-appl}

Before proving the results oulined in Section \ref{sec-setup}, we present 
two applications which show the usefulness of {\bf (GLM2)} in deriving the 
statistical properties of intermittent maps preserving an infinite \me. 

\subsection{Equidistribution of hitting times in residue classes}

Let $T: \ui \into \ui$ be a map satisfying the assumptions of 
Theorem \ref{main-thm-ui}, or the Farey map, or any map of the same
type for which {\bf (GLM2)} holds; cf.\ Remarks \ref{rk-farey} and
\ref{rk-pm2}.

In order to study the intermittent behavior of these maps in
quantitative terms, one looks at how much time the typical \o\ spends in
a neighborhood of the fixed point. The choice of the neighborhood 
is not important, so one usually picks the Markov interval $I_0$. Thus, an 
\ob\ of interest is the \emph{hitting time} of a point $x \in \ui$ to $J := 
\ui \setminus I_0 = [a_1, 1)$:
\begin{equation}
  H(x) := \min \rset{ k\ge 0} {T^k(x) \in J}.
\end{equation}
It is clear that, with the exception of countably many points $x$,
$H(T^n(x))$ is well-defined for all $n \in \N$. We denote by $\rat$
the full-\me\ subset of $\ui$ where the $H \circ T^n$ is well-defined 
for all $n$. 

Consider the level sets of $H$, i.e., $B_k := \rset{x\in \rat}  {H(x)=k}$, 
with $k \in \N$. They form a partition of $\ui$ (mod $\leb$) such that
$B_0 = J \cap \rat$ and, for $k\ge 1$, $B_k \subset I_0$. Also
for $k\ge 1$, $T|_{B_k}$ is a diffeomorphism $B_k \into B_{k-1}$. 
Now, take $x \in \rat$ and consider its \emph{itinerary} 
$( \ell_n ) = ( \ell_n(x) )_{n \in \N}$ w.r.t.\ the partition $\{ B_k \}_{k \in \N}$.
This means that $T^n(x) \in B_{\ell_n}$, for all $n \in \N$. The expansivity 
of $T$ implies that the mapping $x \mapsto (\ell_n)$ is injective, that is, 
equal itineraries correspond to equal points in $\ui$. 

\begin{remark}
  In the case of the Farey map, $\rat = \ui \setminus \Q$ and 
  \begin{displaymath}
    B_k = \left( \frac1{k+2} , \frac1{k+1} \right) \setminus \Q. 
  \end{displaymath}
  The sets $C_k = B_{k-1}$ 
  ($k \ge 1$) are sometimes called Farey cylinders. The itinerary $(\ell_n)$ 
  of a point $x \in \ui \setminus \Q$ is related to its continued fraction 
  expansion $[a_1, a_2, a_3, \ldots]$ as follows:
  \begin{displaymath}
    (\ell_0, \ell_1, \ell_2, \ldots) =
    (a_1-1, a_1-2, \ldots, 0, a_2-1, a_2-2, \ldots, 0, a_3-1, \ldots).
  \end{displaymath}
  Notice that, since $x$ is irrational, the expansion $[a_1, a_2, a_3, 
  \ldots]$ is infinite. (We ask the reader to forgive the abuse of notation 
  whereby in the confines of this remark $a_j$ denotes a digit in the 
  continued fraction expansion, while in the rest of the paper it denotes 
  a point in $[0,1]$.)
\end{remark}

Coming back to the general case, one can use the partition $\{ B_k \}$ 
to construct global \ob s. Given $q \in \Z^+$ and $f_j \in \C$, for 
$0 \le j \le q-1$, denote by $F: \ui \into \C$ the step \fn\ defined 
($\leb$-almost everywhere) by the relation:
\begin{equation} \label{step}
  F (x) = f_j \quad \Longleftrightarrow \quad x \in B_k,\: \mbox{ with } k 
  \equiv j \: (\mbox{mod } q).
\end{equation}

\begin{proposition} \label{prop-hitting}
  Any $F: \ui \into\C$ defined as in {\rm (\ref{step})} is a global \ob\ with 
  \begin{displaymath}
    \avg(F)= \frac1q \sum_{j=0}^{q-1} f_j .
  \end{displaymath}
\end{proposition}

\proof In Section \ref{subs-pf-hitting} of Appendix \ref{sec-tech}.

\skippar

An example of interest, given the discussion at the beginning of the
section, is the global \ob\ $H_q: \ui \into \{0, 1, \ldots, q-1 \}$ given by 
$f_j = j$. The previous proposition shows that  $\avg(H_q)
= (q-1)/2$.
Observe that, for all $x \in \rat$ and $n\in \N$, $H_q(T^n(x)) = 
\ell_n(x) \: (\mbox{mod } q)$.

We want to study the limiting distribution of $H_q \circ T^n$, seen as a
random variable of $x \in \rat$. For this we must specify a \pr\ on $\rat$. 
The invariant \me\ $\mu$ itself is not an option because it is infinite. 
However, since $\mu$ is the reference \me\ of the \dsy, it is reasonable
to use the \pr\ \me\ $\mu_g$ defined by a certain \emph{density}
$g$ relative to $\mu$. In other words, given $g \in L^1(\ui, \mu)$, with 
$g \ge 0$ and $\mu(g)=1$, we consider the \me\ $\mu_g$ such that 
$d\mu_g / d\mu = g$. By Theorem \ref{thm-thaler}\emph{(a)}, $\mu_g$ 
is absolutely continuous w.r.t.\ the Lebesgue \me\ on $\ui$, so it makes 
no difference to think of it as a \me\ on $\ui$ or $\rat$.

It would be desirable for the limiting distribution not to depend on
$g$. We adapt a definition found in \cite[\S3.6]{a}.

\begin{definition}
  Let $F_n$ be a sequence of measurable \fn s $\ui \into \R$, and $X$
  a random variable on some \pr\ space $(\Omega, \P)$. We say that
  $F_n$ converges to $X$ in \emph{strong distributional sense}, as
  $n \to \infty$, if the 
  distribution of $F_n$ w.r.t.\ $\mu_g$ converges to that of $X$, for all 
  densities $g$. In other words, for all \pr\ \me s $\nu \ll \mu$ and 
  all continuous bounded \fn s $\Psi: \R \into \R$,
  \begin{displaymath}
    \lim_{n \to \infty} \int_0^1 (\Psi \circ F_n) \, d\nu = \int_\Omega (\Psi 
    \circ X) \, d\P.
  \end{displaymath}
\end{definition}

\begin{proposition} \label{prop-app1}
  As $n \to \infty$, $H_q \circ T^n$ converges in strong distributional 
  sense to the uniform random variable on the set $\{ 0, 1, \ldots, q-1 \}$.
\end{proposition}

\proof We achieve the result by showing the pointwise convergence
of the corresponding characteristic \fn s. 

The characteristic \fn\ of $H_q \circ T^n$, relative to $\mu_g$,
is given by
\begin{equation}
  \varphi_{n,g}(\t) := \mu_g \!\left( e^{i \t H_q \circ T^n} \right) =
  \mu \! \left( \left( (e^{i \t H_q}) \circ T^n \right) g \right) .
\end{equation}
By Proposition \ref{prop-hitting}, $e^{i \t H_q}$ is a global \ob\ 
with $\avg( e^{i \t H_q} ) = q^{-1} \sum_{j=0}^{q-1} e^{i \t j}$. On the
other hand,
{\bf (GLM2)} implies that, for all densities $g$,
\begin{equation}
  \lim_{n \to \infty} \varphi_{n,g}(\t) = \avg( e^{i \t H_q} ) \mu(g) 
  = \frac1q \sum_{j=0}^{q-1} e^{i \t j}, 
\end{equation}
which is the characteristic \fn\ of the uniform random variable on the 
set $\{ 0, 1, \ldots, q-1 \}$.
\qed

In view of the previous considerations, the above result gives a meaning,
within the scope of infinite \erg\ theory, to the phrase `losing memory
of the initial conditions'. For all choices $\nu \ll \mu$ of the randomness 
of the initial conditions, the $n^\mathrm{th}$ hitting time $\ell_n$, when 
considered mod $q$, converges to the uniform random variable on 
$\{ 0, 1, \ldots, q-1 \}$,  as $n \to \infty$. This is the ``most random'' behavior 
for an \ob\ defined mod $q$.

\subsection{Averaging does not tighten distributions}

The next application is very general and applies to all maps for which 
we have established {\bf (GLM2)} and to a wide class of global \ob s. 

\begin{definition} \label{def-mu-dist}
  For $x,y \in \ui$, respectively $\rp$, let $[x,y]$ denote the closed 
  interval of endpoints $x$ and $y$, irrespective of their order. If $\mu$ is 
  a Lebesgue-equivalent \me\ in $\ui$, respectively $\rp$, the expression
  \begin{displaymath}
    d_\mu(x,y) := \mu( [x,y] )
  \end{displaymath}
  defines a distance in $\ui$, respectively $\rp$, which we call the 
  \emph{$\mu$-distance}.
\end{definition}

Observe that $d_\leb$ is the standard Euclidean distance. In the rest of the 
paper we write that a \fn\ is $d_\mu$-uniformly continuous if it is uniformly 
continuous w.r.t.\ $d_\mu$.

\begin{proposition} \label{prop-app2}
  Let $T$ be a map $\ui \into \ui$ satisfying {\rm (A1)-(A8)} or a map 
  $\rp \into \rp$ satisfying {\rm (B1)-(B6)}, with $\mu$ denoting the
  invariant \me\ (in the latter case, $\mu = \leb$). Let $F$ be a 
  $d_\mu$-uniformly continuous global \ob, taking values in $\R$, such 
  that the infinite-volume average $\avg( e^{i \t F} )$ exists for all 
  $\t \in \R$. Then:
  \begin{itemize}
  \item[(a)] As $n \to \infty$, $F \circ T^n$ converges in strong distributional 
    sense to the random variable $X$ whose characteristic \fn\ is 
    $\varphi_X(\t) := \avg( e^{i \t F} )$;
    
  \item[(b)] For $k \in \Z^+$, denote by 
    \begin{displaymath}
      \mathcal{A}_k F := \frac1{k} \sum_{j=0}^{k-1} F \circ T^j 
    \end{displaymath}
    the $k^{th}$ partial Birkhoff average of $F$. For any fixed $k$,
    as $n \to \infty$, $\mathcal{A}_k F \circ T^n$ converges in strong 
    distributional sense to the same random variable $X$ defined in part 
    (a);
    
  \item[(c)] There exists a diverging sequence $(k_n) \subset \Z^+$ 
    such that $\mathcal{A}_{k_n} F \circ T^n$ converges in strong 
    distributional sense to the variable $X$.  
  \end{itemize}
\end{proposition}

\proof Before starting the proof, we remark that here we have restricted 
to real-valued global \ob s for mere reasons of simplicity. The 
proposition can be easily extended to complex-valued \ob s with the 
suitable modifications.

Statement \emph{(a)} is shown exactly as in the proof of
Proposition \ref{prop-app1} with $F$ in lieu of $H_q$, using that 
$\varphi_X(\t) = \avg( e^{i \t F} )$ exists by hypothesis.

\begin{remark}
  Notice that \emph{(a)} follows directly from {\bf (GLM2)} with
  $e^{i \t F}$ in place of $F$: the hypotheses that $F$ is itself a 
  global \ob\ and that it is $d_\mu$-uniformly continuous are not 
  needed here. More importantly, the argument applies to all types 
  of maps.
\end{remark}

For part \emph{(b)} we need the following lemma, whose proof we
present at the end of Section \ref{sec-pfs-easy}.

\begin{lemma} \label{lem-unif-cont}
  Let $F \in \go$ be $d_\mu$-uniformly continuous and $\Theta: 
  \C^k \into \C$ continuous, for some $k \in \Z^+$. If 
  $\avg( \Theta(F, \ldots, F) )$ exists, then $\avg( \Theta(F \circ T^{n_1}, 
  \ldots, F \circ T^{n_k}) )$ exists for all $n_1, \ldots, n_k \in \N$ and it
  equals $\avg( \Theta(F, \ldots, F) )$.
\end{lemma}

We apply the lemma with 
$\Theta(z_1, \ldots, z_k) := e^{i\t (z_1 + \cdots + z_k) / k}$ and 
$n_j = j-1$. This shows that $\avg( e^{i \t \mathcal{A}_k F } )$ exists 
and equals $\avg( e^{i \t F} ) = \varphi_X(\t)$. Then statement 
\emph{(b)} follows directly from \emph{(a)}. 

As for assertion \emph{(c)}, fix a density $g$ and a positive integer
$k$. Part \emph{(b)} guarantees that there exists a natural number
$\bar{n}_k$ such that
\begin{equation} \label{ult10}
  \left| \mu_g \! \left( e^{i \t \mathcal{A}_k F \circ T^n } \right) -
  \varphi_X(\t) \right| \le 2^{-k}
\end{equation}
for all $n \ge \bar{n}_k$ and all $\t \in E_k := \{ -k, -k+2^{-k}, \ldots, 
k-2^{-k}, k \}$. We can always assume that $\bar{n}_k \nearrow \infty$.
Let $(k_n)_{n \in \N}$ be the following generalized inverse of 
$(\bar{n}_k)_{k \in \Z^+}$:
\begin{equation}
  k_n := \max \rset{1 \le k \le n} {\bar{n}_k \le n}.
\end{equation}
By construction, $n \ge \bar{n}_{k_n}$ for all $n \ge 0$. This fact
and (\ref{ult10}) imply that, for all $\t \in \bigcup_k E_k$, i.e., for
all dyadic rationals $\t$,
\begin{equation} \label{ult30}
  \lim_{n \to \infty} \, \mu_g \! \left( e^{i \t \mathcal{A}_{k_n} F 
  \circ T^n } \right) = \varphi_X(\t).
\end{equation}

The limit is easily extended to all $\t \in \R$, because $F \in L^\infty$
and the random variables $\mathcal{A}_{k_n} F \circ T^n$ are tight. 
A direct proof of this claim is easy, so we give it for the sake of 
completeness. For $\t \in \R$ and $j \in \Z^+$, with $j$ sufficiently 
large, let $\bar{\t}_j$ be an element of $E_j$ that achieves the minimum 
distance from $\t$. Thus $| \t - \bar{\t}_j | \le 2^{-j-1}$. It follows that
\begin{align}
  \left| \mu_g \! \left( e^{i \t \mathcal{A}_{k_n} F \circ T^n } \right) -
  \mu_g \! \left( e^{i \bar{\t}_j \mathcal{A}_{k_n} F \circ T^n } \right) 
  \right| &\le \mu_g \! \left( \left| e^{i (\t - \bar{\t}_j) 
  \mathcal{A}_{k_n} F \circ T^n } - 1 \right| \right) \nonumber \\
  & \le \mu_g \! \left( \left| ( \t - \bar{\t}_j ) \,
  \mathcal{A}_{k_n} F \circ T^n \right| \right) \label{ult40} \\[2pt]
  &\le 2^{-j-1} \, \| F \|_\infty . \nonumber 
\end{align}
Given $\eps>0$, choose $j$ so large that $2^{-j-1} \,
\| F \|_\infty \le \eps/3$ and
\begin{equation} \label{ult50}
  \left| \varphi_X( \bar{\t}_j ) - \varphi_X(\t) \right| \le \frac{\eps}3.
\end{equation}
The first condition implies that the rightmost term of (\ref{ult40}) does 
not exceed $\eps/3$ \emph{for all} $n$. The second condition is 
possible because of the continuity of the characteristic \fn. Now apply 
(\ref{ult10}) with $\t := \bar{\t}_j$ and $k := k_n$: its l.h.s.\ can be 
made smaller than or equal to $\eps/3$ for all sufficiently large $n$. 

Combining all these inequalities proves (\ref{ult30}) for an arbitrary 
density $g$, ending the proof of part \emph{(c)}. 
\qed

Statements \emph{(b)} and \emph{(c)} of Proposition \ref{prop-app2}
are in sharp contrast to what happens in \m\ systems preserving a 
probability \me\ $\mu$. In \emph{all} such cases, consider a non-constant
bounded \fn\ $f$ and denote by $X$ the random variable given by $f$
w.r.t.\ the probability $\mu$, in other words, the one determined by the
characteristic \fn\ $\varphi_X(\t) := \mu( e^{i \t f} )$. We have:
\begin{enumerate}
  \item As $n \to \infty$, $\mathcal{A}_k f \circ T^n$ converges in strong 
  distributional sense to a variable that, for large $k$, has a smaller 
  variance than $X$. 
  
  \item For any diverging sequence $(k_n) \subset \Z^+$, 
  $\mathcal{A}_{k_n} f \circ T^n$ does not converge in strong distributional 
  sense to $X$.
  
  \item There exists a diverging sequence $(k_n) \subset \Z^+$
  such that $\mathcal{A}_{k_n} f \circ T^n$ converges in strong 
  distributional sense to the constant $\mu(f)$.
\end{enumerate}

These claims are easily proved. In fact, for any density $g$ 
and any $\t \in \R$, using \m, we have
\begin{equation} \label{ult55}
  \lim_{n \to \infty} \mu_g \! \left( e^{i \t \mathcal{A}_k f \circ T^n} 
  \right) := \lim_{n \to \infty} \mu \! \left( e^{i \t \mathcal{A}_k 
  f \circ T^n} g \right) = \mu \! \left( e^{i \t \mathcal{A}_k f} \right),
\end{equation}
so the limiting variable in statement 1 is given by the \fn\ 
$\mathcal{A}_k f$ w.r.t.\ the \pr\ $\mu$. On the other hand, again by
\m, 
\begin{equation}
  \left| \mu\! \left( [f \circ T^j - \mu(f)] \, [f - \mu(f)] \right) \right|
  < \mu\! \left( [f - \mu(f)]^2 \right),
\end{equation}
for all sufficiently large $j$ (observe that the above r.h.s.\ is positive because
$f$ is non-constant). This and the invariance of $\mu$ imply that, for $k$ large 
enough, 
\begin{equation}
  \mu \! \left( [ \mathcal{A}_k f - \mu(f) ]^2 \right) < \mu \! \left( [f - \mu(f) ]^2 \right) , 
\end{equation}
giving our first claim.

For the second claim let us chose the density 1; in other words, let us consider 
$\mathcal{A}_{k_n} f \circ T^n$ as a random variable w.r.t.\ $\mu$. 
Since $\mu$ is invariant, the distribution of  $\mathcal{A}_{k_n} f \circ T^n$ 
is the same as that of $\mathcal{A}_{k_n} f$. By \erg ity, the latter variable
converges almost everywhere, and thus in distribution, to the constant
$\mu(f)$, which cannot be equal to the non-constant variable $X$.

For the third claim we proceed as in the proof of \emph{(c)}. Using
(\ref{ult55}) we find a suitable sequence $(k_n)$ such that, for all dyadic 
rationals $\t$,
\begin{equation} \label{ult100}
  \lim_{n \to \infty} \left| \mu_g \! \left( e^{i \t \mathcal{A}_{k_n} f 
  \circ T^n } \right) - \mu \! \left( e^{i \t \mathcal{A}_{k_n} f}
  \right) \right| = 0.
\end{equation}
The limit is then extended to all $\t \in \R$ by tightness, as shown 
earlier. On the other hand, by \erg ity, $\mathcal{A}_{k_n} f \to
\mu(f)$ $\mu$-almost everywhere, as $n \to \infty$, implying
that 
\begin{equation} \label{ult110}
  \lim_{n \to \infty} \mu \! \left( e^{i \t \mathcal{A}_{k_n} f}
  \right) = e^{i \t \mu(f)}.
\end{equation}
Statement 2 then follows from (\ref{ult100}) and (\ref{ult110}). 

As a final comment, Proposition \ref{prop-app2} is a consequence of the 
fact that any absolutely continuous finite \me\ is eventually pushed to a 
neighborhood of the fixed point. This only occurs when the fixed point is 
strongly neutral, giving rise to an infinite invariant \me.

\section{First proofs}
\label{sec-pfs-easy}

The rest of the paper is largely devoted to the proofs of the results 
presented in Section \ref{sec-setup}. In this section we deal with the 
simpler results, Propositions \ref{prop-inv-ui}, \ref{prop-other-ui}, 
\ref{prop-inv-rp} and \ref{prop-other-rp}. In fact we will only write the proofs 
of the first two, as the other two are analogous---indeed easier, as they 
involve the Lebesgue \me\ instead of $\mu$. At the end of the section
we also give the proof of Lemma \ref{lem-unif-cont}, which was left behind.

\skippar

\proofof{Proposition \ref{prop-inv-ui}} The proposition will be proved
once we show that, for all $n \in \N$,
\begin{equation} \label{inv10}
  \lim_{\ivui} \frac{ \mu(T^{-n} V \symmdiff V) } {\mu(V)} = 0,
\end{equation}
where $\symmdiff$ denotes the symmetric difference of two sets. 
In fact, the invariance of $\mu$ and the boundedness of $F$ imply that
\begin{equation}
  \frac1 {\mu(V)} \int_V F \, d\mu = \frac1 {\mu(V)} \int_{T^{-n} V} 
  (F \circ T^n) \, d\mu = \frac1 {\mu(V)} \int_V (F \circ T^n) \, 
  d\mu + \epsilon(V),
\end{equation}
where $\epsilon(V)$ is an error term that is bounded above by 
$\|F\|_\infty \, \mu(T^{-n} V \symmdiff V) / \mu(V)$.

So it remains to verify (\ref{inv10}) in our specific case. Since $T^n$ is 
again a piecewise smooth Markov map with countably many surjective 
branches and an indifferent fixed point at 0, we can assume $n=1$. 

Write $V := [a,1)$. The infinite-volume limit is $a \to 0^+$. Using (A1) 
and (A4) we have
\begin{equation} \label{inv20}
  T^{-1} V = \bigcup_{j \in \ind} \br_j^{-1} V = [\br_0^{-1}(a), a_1) \cup
  \bigcup_{j \ge 1} \br_j^{-1} V.
\end{equation}
Observe that $[\br_0^{-1}(a), a_1) \supset [a, a_1)$ and $\br_j^{-1} V
\subset I_j$. Thus,
\begin{equation} \label{inv30}
  T^{-1} V \symmdiff V = [\br_0^{-1}(a), a) \cup \bigcup_{j \ge 1} 
  ( I_j \setminus \br_j^{-1} V ).
\end{equation}
The relation $\mu(T^{-1} V) = \mu(V)$ implies that
\begin{equation} \label{inv40}
  \mu( [\br_0^{-1}(a), a) ) = \sum_{j\ge 1} \mu(I_j \setminus \br_j^{-1} V).
\end{equation}
Observe that $\mu$ is a finite \me, when restricted to $\bigcup_{j\ge 1} I_i$, 
and each $I_j \setminus \br_j^{-1} V$ decreases to the empty set, as $a$ 
decreases to 0. Therefore (\ref{inv40}) vanishes for $a \to 
0^+$. Applied to (\ref{inv30}), this shows that $\mu(T^{-1} V \symmdiff V) 
\to 0$, as $\ivui$, implying (\ref{inv10}).
\qed

\proofof{Proposition \ref{prop-other-ui}} {\bf (GLM1)} and 
{\bf (LLM)} come from the exactness of $T$ and 
\cite[Thm.\ 3.5]{lpmu}.

In order to show that no form of global-global \m\ holds, let us pick a 
real-valued, $d_\mu$-uniformly continuous global \ob\ $F$ (cf.\ 
Definition \ref{def-mu-dist}) such that $\avg(F^2)$ exists and is 
different from  $(\avg(F))^2$. One 
example is $F(x) := \sin(\Phi(x))$, where $\Phi$ is the \fn\ defined in 
Section \ref{subs-rp}, mapping $(\ui, \mu)$ onto $(\rp, \leb)$: one can 
easily verify that $\avg(F) = \avgleb(\sin) = 0$ and $\avg(F^2) = 
\avgleb(\sin^2) = 1/2$.

To this \ob\ we apply Lemma \ref{lem-unif-cont}, which we stated
in Section \ref{sec-appl} and will prove momentarily. (The proof will
not involve any of the results of Section \ref{sec-appl}, so there is no
circular reasoning.) Specifically we apply the lemma with $k := 2$, 
$\Theta(z_1, z_2) := z_1 z_2$, $n_1 := n$ and $n_2 := 0$. Thus, 
$\exists \, \avg ( (F \circ T^n) F) = \avg(F^2) \ne (\avg(F))^2$. 
This contradicts both {\bf (GGM1)} 
and {\bf (GGM2)}.

Finally, {\bf (GLM3)} does not hold because otherwise Proposition 2.4 
of \cite{lpmu} (whose hypotheses hold here) would imply {\bf (GGM2)}.
This concludes the proof of Proposition \ref{prop-other-ui}.
\qed

\skippar

\proofof{Lemma \ref{lem-unif-cont}} Once again, we only write the proof 
for the case of maps $T: \ui \into \ui$ satisfying (A1)-(A4). The case 
$T: \rp \into \rp$ satisfying (B1)-(B5) is completely analogous (using $\leb$
in place of $\mu$).

Since $\Theta$ is continuous,
it is uniformly continuous on every compact set of $\C^k$. In 
particular, for all $\eps >0$, there exists $\delta > 0$ such that,
every time $|z_j|, |w_j| \le \|F\|_\infty$ and $|z_j - w_j| \le \delta$
(for $1 \le j \le k$), one has
\begin{equation} \label{uc2}
  \left| \Theta(z_1, \ldots, z_k) - \Theta(w_1, \ldots, w_k) \right| \le
  \eps.
\end{equation}
By the uniform continuity of $F$, we can find $\gamma > 0$ such 
that 
\begin{equation} \label{uc4}
  d_\mu(x,y) \le \gamma \quad \Longrightarrow \quad 
  | F(x) - F(y) | \le \delta.
\end{equation}

Now we claim that, for any $n \in \N$ and $\gamma > 0$, there exists 
$a' \in (0,1)$ such that, for all $x \in (0,a']$, $d_\mu (x, T^n(x)) \le \gamma$. 
To establish this claim, note that we can suppose without loss of 
generality that $n=1$ and use arguments from the proof of 
Proposition \ref{prop-inv-ui}. So, for $a \in (0, a_1]$, set
$V' := [T(a),1)$ and proceed as in 
(\ref{inv20})-(\ref{inv40}), with $T(a)$ in lieu of $a$. Since $a \in 
(0, a_1]$, we have that $\br_0^{-1}(T(a)) = a$, whence
$d_\mu(a, T(a)) := \mu( [a,T(a)] ) \searrow 0$, as $a \searrow 0$. 
Finally, let  $a' \in (0, a_1]$ be uniquely defined by 
$d_\mu(a', T(a')) = \gamma$. By the monotonicity of the limit
in $a$, $d_\mu(x, T(x)) \le \gamma$, for all $x \in (0,a']$.

We make a repeated use of the above claim with $n = n_j$,
for $j = 1, 2, \ldots, k$. In each case, we obtain some $a'_j$ such
that $d_\mu (x, T^{n_j}(x)) \le \gamma$, for all $x \in (0, a'_j]$. Set 
$\bar{a} := \min_{1 \le j \le k} \{ a'_j \}$. In view of (\ref{uc2})-(\ref{uc4}), 
we see that, for $x \in (0, \bar{a}]$,
\begin{equation} \label{uc10}
  \big|  \Theta( F(x), \ldots, F(x) ) -
  \Theta( F(T^{n_1}(x)), \ldots, F(T^{n_k}(x)) ) \big| \le \eps.
\end{equation}

Recall the notation (\ref{mu-v-ui})-(\ref{def-avg-ui}). For $a < \bar{a}$, 
\begin{equation} \label{uc20}
\begin{split}
  & \mu_{[a,1)} ( \Theta (F \circ T^{n_1}, \ldots, F \circ T^{n_k}) ) \\[4pt]
  & \qquad = \frac1 {\mu( [a,1) )} \int_a^{\bar{a}} \Theta 
  (F \circ T^{n_1}, \ldots, F \circ T^{n_k}) \, d\mu \\
  & \qquad \qquad + \frac1 {\mu( [a,1) )} \int_{\bar{a}}^1  
  \Theta (F \circ T^{n_1}, \ldots, F \circ T^{n_k}) \, d\mu.
\end{split}
\end{equation}
As $a \to 0^+$, the second term of the above r.h.s.\ vanishes. 
Furthermore, by (\ref{uc10}),
\begin{equation} \label{uc30}
  \left| \frac1 {\mu( [a,1) )} \int_a^{\bar{a}} \Theta (F \circ T^{n_1}, \ldots, 
  F \circ T^{n_k}) \, d\mu - \frac1 {\mu( [a,1) )} \int_a^{\bar{a}} \Theta 
  (F, \ldots, F) \, d\mu \right| \le \eps.
\end{equation} 
In analogy with (\ref{uc20}),
\begin{equation} \label{uc40}
\begin{split}
  & \left| \mu_{[a,1)} ( \Theta (F, \ldots, F) ) - \frac1 {\mu( [a,1) )} 
  \int_a^{\bar{a}} \Theta (F, \ldots, F) \, d\mu \right| \\
  & \qquad \le \frac1 {\mu( [a,1) )} \int_{\bar{a}}^1 | \Theta (F, \ldots, F) | 
  \, d\mu,
\end{split}
\end{equation} 
which vanishes as $a \to 0^+$. Putting everything together, we
obtain
\begin{equation} \label{uc50}
  \limsup_{a \to 0^+} \left| \mu_{[a,1)} ( \Theta (F \circ T^{n_1}, \ldots, 
  F \circ T^{n_k}) ) - \mu_{[a,1)} ( \Theta (F, \ldots, F) ) \right| \le \eps.
\end{equation} 
Since $\eps$ is arbitrary, the above limit proves the assertion of Lemma 
\ref{lem-unif-cont}.
\qed

\section{Proof of (GLM2)} 
\label{sec-glm2}

The proof of {\bf (GLM2)} follows the same strategy for both maps on 
$\ui$ and $\rp$. It hinges on the exactness of the maps and the 
existence of a local \ob\ with a certain monotonicity property, see 
Definition \ref{def-per-mon} below. What changes in the two cases is the 
assumptions that are needed to guarantee the existence of this special 
\ob. We will deal with this in Section \ref{sec-pers-mon}.

For the rest of the paper we use the bracket notation to indicate the 
integral product of a global \ob\ and a local \ob, w.r.t.\ to the invariant 
\me. More precisely, for $F \in L^\infty$ and $g \in L^1$, we define
\begin{equation} \label{coupling-ui}
  \langle F, g \rangle := \int_0^1 F g \, d\mu,
\end{equation}
if we are working with the space $\ui$, and
\begin{equation} \label{coupling-rp}
  \langle F, g \rangle :=  \int_0^\infty F g \, d\leb,
\end{equation}
if we are working in $\rp$. 

Denote by $P = P_T$ the transfer operator of $T$, relative to the
above coupling. This is defined by the identity $\langle F \circ T,
g \rangle = \langle F, Pg \rangle$. The functional forms of $P$ in
the two cases are given, respectively, in (\ref{trans-op-ui}) and
(\ref{trans-op-rp}).

\begin{definition} \label{def-per-mon}
  We say that the local \ob\ $g$ is \emph{persistently monotonic}
  if, for all $n \in \N$, $P^n g(x)$ is a positive, monotonic \fn\ of $x$. 
\end{definition}

For maps $T: \ui \into \ui$, the above condition reads: $P^n g$ is an 
increasing \fn\ of $\ui$. For maps $T: \rp \into \rp$, it reads: $P^n g$ is 
a decreasing \fn\ of $\rp$. 

The following theorem contains the main idea of the paper.

\begin{theorem} \label{thm-main}
  Let $T$ be a map $\ui \into \ui$ verifying {\rm (A1)-(A4)}, or a map
  $\rp \into \rp$ verifying {\rm (B1)-(B5)}. If $T$ admits a persistently 
  monotonic local \ob, then it is {\bf (GLM2)}.
\end{theorem}

\proof Once again, we prove the result only for $T: \ui \into \ui$, the 
other case being analogous and simpler. 
We use \cite[Lem.~3.6]{lpmu}, which we restate here in a convenient
form.

\begin{lemma} \label{lem-main-pmu}
  Assume that $T$ is exact and $F \in \go$.  If the limit
  \begin{displaymath}
    \lim_{n \to \infty} \mu( (F \circ T^n) g) = \avg(F) \mu(g)  
  \end{displaymath}
  holds for some $g \in \lo$, with $\mu(g) \ne 0$, then it holds for
  all $g \in \lo$.
\end{lemma}

Thus, recalling that $T$ is exact by Theorem \ref{thm-thaler}\emph{(b)},
it suffices to verify the above limit when $g$ is the
persistently monotonic \ob\ provided by the hypotheses
of the theorem. Notice that $\mu(g) = \| g \|_1 > 0$. Without loss of 
generality, we can assume $\| g \|_1 = 1$, otherwise one considers 
$g_1 := g / \| g \|_1$.

It all reduces to prove that, for all $F \in \go$ with $\avg(F)=0$,
\begin{equation} \label{main-goal}
  \lim_{n \to \infty} \langle F, P^n g \rangle = 0.
\end{equation}
In fact, for $\avg(F) \ne 0$, one applies (\ref{main-goal}) to $F_1 := 
F - \avg(F)$, which satisfies $\avg(F_1) = 0$. It follows that
\begin{equation} 
  \lim_{n \to \infty} \langle F \circ T^n, g \rangle = \lim_{n \to \infty}
  \langle F, P^n g \rangle = \avg(F) \mu(g).
\end{equation}
Here one uses that, for $g>0$, $\langle 1, P^n g \rangle = \| P^n g \|_1 
=  \| g \|_1 = \mu(g)$.

So, fix $F \in \go$, with $\avg(F) = 0$, and $\eps > 0$. By definition, cf.\
(\ref{mu-v-ui})-(\ref{def-avg-ui}), there exists $\delta > 0$ such that
\begin{equation} \label{main20}
  \forall a \le \delta, \quad \frac1{\mu( [a,1) )} \left| \int_a^1 F \, d\mu \right| 
  < \frac{\eps}2.
\end{equation}
For $x \in \ui$ and $n \in \N$, set
\begin{equation} \label{main30}
  \gamma_n(x) = \gamma_{n,\delta}(x) := \min\{ P^n g(\delta),
  P^n g(x) \}.
\end{equation}
Since $g$ is persistently monotonic, $\gamma_n$ is a positive, 
increasing \fn, with a plateau on $[\delta,1)$. It is a local
\ob\ because $\| \gamma_n \|_1 \le \| P^n g \|_1 =  \| g \|_1 = 1$. We 
have
\begin{equation} \label{main-goal2}
  \langle F, P^n g \rangle = \int_0^1 F \gamma_n \, d\mu +
  \int_\delta^1 F (P^n g - \gamma_n) \, d\mu =: \mathcal{I}_1 + 
  \mathcal{I}_2.
\end{equation}

To estimate $\mathcal{I}_2$, let us notice that
\begin{equation} 
  0 \le \int_\delta^1 (P^n g - \gamma_n) \, d\mu \le \int_\delta^1 P^n g 
  \, d\mu = \langle 1_{[\delta,1)} , P^n g \rangle.
\end{equation}
Since the \sy\ is {\bf (LLM)}  (Proposition \ref{prop-other-ui}), the
rightmost term above vanishes, as $n \to \infty$. Thus, for all sufficiently 
large $n$,
\begin{equation} \label{main50}
  | \mathcal{I}_2 |  \le
  \| F \|_\infty  \int_\delta^1 (P^n g - \gamma_n) \, d\mu \le \frac{\eps}2.
\end{equation}

Let us consider $\mathcal{I}_1$. For $0 \le r < \gamma_n(\delta) =
P^n g (\delta)$, the expression
\begin{equation}
  \gamma_n^{-1}(r) := \inf \rset{x \in \ui} {\gamma_n(x) \ge r}
\end{equation} 
defines the generalized inverse of $\gamma_n$, which is an
increasing \fn\ of $r$. Using a trick and Fubini's Theorem, we can write
\begin{equation}
  \mathcal{I}_1 = \int_0^1 F(x) \left( \int_0^{\gamma_n(x)} \!\!\ dr 
  \right) \mu(dx) = \int_0^{\gamma_n(\delta)} \left( 
  \int_{\gamma_n^{-1}(r)}^1 F(x) \, \mu(dx) \right) dr.
\end{equation} 
Therefore, using (\ref{main20}) with $a := \gamma_n^{-1}(r)$ and 
observing that $\gamma_n^{-1}(r) \le \delta$ by construction, we get
\begin{equation} \label{main80}
\begin{split}
  | \mathcal{I}_1 | &\le \int_0^{\gamma_n(\delta)} \left| 
  \int_{\gamma_n^{-1}(r)}^1 F(x) \, \mu(dx) \right| dr \\
  & \le \frac{\eps}2 \int_0^{\gamma_n(\delta)} 
  \int_{\gamma_n^{-1}(r)}^1 \mu(dx) \ dr \\
  &= \frac{\eps}2 \int_0^1 \int_0^{\gamma_n(x)} \!\!\ dr \
  \mu(dx) \\
  &= \frac{\eps}2 \mu(\gamma_n) \le \frac{\eps}2.
\end{split}
\end{equation} 
The trick that we have used effectively consists in disintegrating 
the density $\gamma_n$ in infinitely many horizontal slices, one for 
each value of $r$. Each slice corresponds to an infinitesimal
multiple of the \pr\ distribution $\mu_{[\gamma_n^{-1}(r) , 1)}$,
relative to which $F$ has an almost zero average.

The estimate (\ref{main80}) holds uniformly in $n$. Together with 
(\ref{main50}) and (\ref{main-goal2}), it proves (\ref{main-goal}).
\qed

One might wonder how the above arguments relate to the 
technique used to prove global-local \m\ for the uniformly expanding
maps of \cite{lmmaps}. In general they do not: the proof of 
{\bf (GLM2)} for quasi-lifts on $\R$ uses a different idea based on 
the invariance of such maps for the action of $\Z$. However, 
for the special example of \cite[Sect.~4.3]{lmmaps}, the author 
employs an argument which is the discrete equivalent of the 
slicing of the density described earlier.

\section{Persistently monotonic local observables}
\label{sec-pers-mon}

In this section we establish the existence of persistently monotonic 
local \ob s in the two cases considered. Together with Theorem 
\ref{thm-main}, this will prove Theorems \ref{main-thm-ui} and 
\ref{main-thm-rp}.

\subsection{Case of the unit interval} 
\label{sub-unitint}

For a map $T$ verifying (A1)-(A4), the transfer operator $P=P_T$ 
relative to the coupling (\ref{coupling-ui}) reads
\begin{equation} \label{trans-op-ui}
  (Pg)(x) = \frac{1}{h(x)}\, \sum_{j\in \ind}\, |\ibr'_j(x)|\, h(\ibr_j(x))\, 
  g(\ibr_j(x)),
\end{equation}
where $h = d\mu / d\leb$. If, with a harmless abuse of notation, we let 
$P$ act on $L^\infty$ too, we see that $P1=1$, which is 
equivalent to the invariance of $\mu$.

\begin{theorem}\label{exis-pm-ui}
  Let $T: \ui \into \ui$ satisfy assumptions {\rm (A1)-(A8)} of Sections
  \ref{subs-ui} and \ref{subs-im-ui}. Then $T$ admits a persistently 
  monotonic local observable.
\end{theorem}

\proof We claim that if $g: \ui \into \R$ is a differentiable, positive, 
increasing, concave local \ob, then the same holds for $Pg$. So, by 
induction, any $g$ with these features is such that $P^n g$ is a
positive and increasing local \ob\ for all $n \in \N$, proving the 
theorem.

We prove the claim by means of the following technical lemma,
whose proof is given in Section \ref{subs-increasing-ui} of
Appendix \ref{sec-tech}.

\begin{lemma} \label{increasing-ui}
Take a differentiable, increasing and concave \fn\ $g: \ui \into \R$.
Let $\ibr_0, \ibr_1: \ui \into \ui$ be twice differentiable and such that
\begin{itemize}
\item[{\rm (H1)}] $\ibr_0' \ge 0$;
\item[{\rm (H2)}] $\max \ibr_0 \le \min \ibr_1$;
\item[{\rm (H3)}] $\ibr_0' + \ibr_1' \ge 0$;
\item[{\rm (H4)}] $\ibr_0'' \le 0$ and $\ibr_0'' + \ibr_1'' \le 0$.
\end{itemize}
Also, let $\aa,\bb : \ui \into \R^+$ be differentiable and such that
\begin{itemize}
\item[{\rm (H5)}] $\aa+\bb=1$;
\item[{\rm (H6)}] $\aa \ge \bb$;
\item[{\rm (H7)}] $\aa$ is decreasing and convex.
\end{itemize}
Then
\begin{equation} \label{f1}
  g_1 := \aa \, (g \circ \ibr_0) + \bb \, (g \circ \ibr_1)
\end{equation}
is differentiable, increasing and concave.
\end{lemma}

It is immediate to verify that (A1)-(A6) imply the hypotheses (H1)-(H4) 
of the lemma. Let us set
\begin{equation}
  \aa := \ibr'_0 \, \frac{h \circ \ibr_0}{h} ; \qquad \bb := -\ibr'_1 \, 
  \frac{h \circ \ibr_1}{h}.
\end{equation}
With these definitions, in view of (\ref{trans-op-ui}) and (\ref{f1}),
and using (A4) and (A5), we have that $g_1 = Pg$. The identity
$P1=1$ gives (H5), while (A7) and (A8) imply, respectively, (H7) 
and (H6). Finally, $\aa$ is differentiable by (A7) and $\bb$ is 
differentiable by (H5). So Lemma \ref{increasing-ui} can be applied.

To finish the proof of the claim it remains to
observe that if $g$ is a positive local \ob, then $Pg$ is also a positive 
local \ob, because $P$ is the transfer operator.
\qed

\subsection{Case of the half line} 
\label{sub-halfline}

For a map $T$ verifying (B1)-(B5), the transfer operator $P=P_T$ relative to the 
coupling (\ref{coupling-rp}) reads
\begin{equation} \label{trans-op-rp}
  (Pg)(x) = \sum_{j\in \ind}\, |\ibr'_j(x)|\, g(\ibr_j(x))
\end{equation}
Once again, $P1=1$. Comparing (\ref{trans-op-rp}) with (\ref{trans-op-ui}), it is clear 
why assumption (B5) simplifies our proof here.

\begin{theorem} \label{exis-pm-rp}
  Let $T: \rp \into \rp$ satisfy assumptions {\rm (B1)-(B6)} of Sections 
  \ref{subs-rp} and \ref{subs-im-rp}.Then $T$ admits a persistently 
  monotonic local observable.
\end{theorem}

\proof As in the proof of Theorem \ref{exis-pm-ui}, we use a recursive
argument. Specifically, we show that if $g: \rp \into \R$ is a differentiable, 
positive, decreasing local \ob, then so is $Pg$.

By (B6) we know that $\ibr'_j > 0$, for all $j \in \ind$, hence, in terms of \fn s
$\rp \into \R$, (\ref{trans-op-rp}) becomes
\begin{equation} \label{trans-op-rp2}
  Pg = \sum_{j\in \ind} \ibr'_j \, ( g \circ \ibr_j ).
\end{equation}
The \fn\ $Pg$ is a positive local \ob\ by general properties of $P$, and is 
differentiable because $\ibr_j$ is $C^2$ by (B1). It remains to show 
that $(Pg)' \le 0$.

The invariance of $\leb$, equivalently, the identity $P1=1$, gives $\sum_{j \in \ind} 
\ibr'_j = 1$, whence
\begin{equation} \label{eq-deri-sec}
  \sum_{j\in \ind}\, \ibr''_j = 0.
\end{equation}
Differentiating (\ref{trans-op-rp2}) gives
\begin{equation} \label{eq-bo0}
  (Pg)' = \sum_{j\in \ind} \ibr''_j \, (g \circ \ibr_j) + \sum_{j\in \ind} (\ibr'_j)^2\, 
  (g' \circ \ibr_j).
\end{equation}
Since $g$ is decreasing,
\begin{equation} \label{eq-bo1}
  \sum_{j\in \ind} (\ibr'_j)^2\, (g' \circ \ibr_j) \le 0.
\end{equation}
By definition, $\ibr_j < \ibr_0$, for all $j\ge 1$, implying 
\begin{equation} \label{eq-bo2}
  g \circ \ibr_j \ge g \circ \ibr_0.
\end{equation}
Finally, (B6) ensures that $\ibr''_j \le 0$, for all $j \ge 1$. This, (\ref{eq-bo2}) and 
(\ref{eq-deri-sec}) give
\begin{equation} \label{eq-bo3}
  \sum_{j\in \ind} \ibr''_j \, (g \circ \ibr_j) \le \left( \sum_{j\in \ind} \ibr''_j \right) 
  (g \circ \ibr_0) = 0.
\end{equation}
The (in)equalities (\ref{eq-bo0}), (\ref{eq-bo1}) and (\ref{eq-bo3}) show that
$(Pg)' \le 0$, ending the proof of Theorem \ref{exis-pm-rp}.
\qed

\appendix

\section{Appendix: Exactness for maps of the half-line}
\label{sec-ex}

In this section we prove a generalization of Theorem \ref{thm-ex-leb}
to the case where $T$ preserves an absolutely continuous
(not necessarily infinite) \me\ $\mu$ on $\R^+$. Specifically,
we replace (B5) of Section \ref{subs-rp} with the weaker 
assumption

\begin{itemize}
\item[(B5')] $T$ preserves an absolutely continuous \me\ $\mu$
  such that $d\mu / d\leb$ is positive and locally integrable on
  $\R^+$.
\end{itemize}

\begin{theorem} \label{thm-ex} 
  Under the assumptions {\rm (B1)-(B4)} and {\rm (B5')}, 
  $T$ is conservative and exact.
\end{theorem}

\proof This proof is based on that of \cite[Thm.~2.1]{lsimple}. 
In the following we outline the flow of the proof, but do not reprove
the statements from \cite{lsimple} that apply verbatim here.
We instead concentrate on the arguments that need modification. 

Set $J := \R^+ \setminus I_0$. Assumption (B4) ensures that,
for any $x \in I_0$, $T^n(x)$ decreases until it lands in $J$,
for some $n$. 
Hence, $J$ is a global cross-section, in the sense
that almost every \o\ of the \sy\ intersects it. Then (B5') 
implies that $\mu(J) < \infty$ and that the Poincar\'e Recurrence 
Theorem can be applied to the map induced by $T$ on $J$,
w.r.t.\ invariant \me\ $\mu$. Therefore the \sy\ is conservative.

For the exactness we apply the Miernowski-Nogueira
criterion \cite{mn}:

\begin{proposition} \label{prop-mn} 
  The non-singular, \erg\ \dsy\ $(X, \sca, \nu, \mathcal{T})$ is 
  exact if and only if, $\forall A \in \sca$ with $\nu(A) > 0$, $\exists 
  n = n(A)$ such that $\nu( \mathcal{T}^{n+1}A \cap \mathcal{T}^n 
  A ) > 0$.
\end{proposition}

(See \cite[Sect.~A.2]{lsimple} for a generalization of the above
criterion to the case of non-\erg\ \sy s.) 

\skippar

We use Proposition \ref{prop-mn} with $\mathcal{T} = T$ and $\nu = 
\leb$; that is, from this point forth we use the Lebesgue \me\ $\leb$.
We have already seen that $J = \bigcup_{j \in \ind \setminus \{0\} } I_j$ 
is a global cross section. Given a positive-\me\ set $A$, we claim that
the forward \o\ of a typical $x_0 \in A$ visits some interval 
$I_{\bar{\jmath}}$, with $\bar{\jmath} \in \ind \setminus \{0\}$,
an infinite number of times. In fact, in the opposite case, 
$(T^n(x_0))_{n \ge 0}$ must eventually leave every $I_j$ for good,
implying that $T^n(x_0) \to 0$. However, by conservativity, 
this can only happen for a null set of points.

The typical $x_0 \in A$ is also a point of density 1 for $A$, relative to the
Lebesgue \me\ $\leb$, namely,
\begin{equation} \label{ex-0}
  \lim_{\eps \to 0^+} \leb ( A \,|\, [x_0-\eps, x_0+\eps] ) :=
  \lim_{\eps \to 0^+} \frac {\leb ( A \cap [x_0-\eps, x_0+\eps] )}
  {\leb ( [x_0-\eps, x_0+\eps] )} = 1.
\end{equation}
Moreover, $I_{\bar{\jmath}}$ is a Markov interval
for $T$, which is uniformly expanding away from $+\infty$. 
Therefore, if $(n_k)_{k \ge 1}$ is the sequence of the hitting times of 
$x_0$ to $I_{\bar{\jmath}}$, $T^{n_k}$ maps a smaller and smaller interval
around $x_0$, where the density of $A$ is higher and higher, 
onto $I_{\bar{\jmath}}$. (The small interval in question is of the form
(\ref{ex-5}), see later; cf.\ also (\ref{ex-10}).)

\emph{If} the map has bounded distortion, the above implies that
the density of $T^{n_k} A$ within $I_{\bar{\jmath}}$ also gets higher 
and higher, as $k$ grows. More precisely, in view of (\ref{ex-0}),
\begin{equation} \label{ex-1}
  \lim_{k \to \infty} \leb( T^{n_k} A \,|\, I_{\bar{\jmath}} ) = 1. 
\end{equation}
It follows that $\exists k \in \Z^+$ such that
\begin{equation} \label{ex-2}
  \leb(T^{n_k + 1} A \cap T^{n_k} A) > 0.
\end{equation}
Let us show this. Denote by $B := \br_{\bar{\jmath}}^{-1} I_{\bar{\jmath}}$ 
the preimage of $I_{\bar{\jmath}}$ via the $\bar{\jmath}^\mathrm{th}$
branch of $T$, which is surjective. Then $B$ is a positive-\me\
subset of $I_{\bar{\jmath}}$. By (\ref{ex-1}), for all sufficiently large $k$,
\begin{align}
\label{ex-3}
  \leb( T^{n_k} A \,|\, I_{\bar{\jmath}} ) &> \frac12; \\
\label{ex-4}
  \leb( T^{n_k} A \,|\, B ) &> 1 - \frac1{2D}, 
\end{align}
where $D$ is the distortion coefficient of $T$ (cf.\ Lemma \ref{lem-dist}
below). Applying $T$ to (\ref{ex-4}) gives 
$\leb( T^{n_k+1} A \,|\, I_{\bar{\jmath}} ) > 1/2$, which, together with
(\ref{ex-3}), yields (\ref{ex-2}).

We have thus verified the main hypothesis of Proposition \ref{prop-mn}. 
The proposition also requires that $T$ be \erg. But this is easy to verify: 
if $A$ is a positive-\me\ \emph{invariant} set, (\ref{ex-1}) reads
$\leb ( A \,|\, I_{\bar{\jmath}} ) = 1$, or $A = I_{\bar{\jmath}}$ mod
$\leb$, whence $A = TA = T I_{\bar{\jmath}} = \R^+$ mod $\leb$. 

Therefore, up to details which can be checked in the proof of 
\cite[Thm.~2.1]{lsimple}, it remains to show that $T$ has bounded 
distortion. This part too follows the same line of reasoning as the 
aforementioned proof, although, understandably, some of the 
computations are different. In order to state the needed result we need 
some preparatory material. 

Set $b_0 := a_1$. For $k \ge 1$, let $b_k$ be uniquely defined by
$b_k > b_{k-1}$ and $T(b_k) = b_{k-1}$. Set $I_{-k} := 
(b_{k-1}, b_k)$; evidently, $\pa_- := \{ I_j \}_{j \in \Z^-}$ is a partition 
of $I_0$ (mod $\leb$). So $\pa_o := \pa_- \cup \pa \setminus \{ I_0 \}$ 
is a partition of $\R^+$, with index set  $\ind_o := \Z^- \cup \ind 
\setminus \{0\}$. 
$\pa_o$ is a refinement of $\pa$, and still a Markov partition for 
$T$, because $T( I_{-1} ) = J = \bigcup_{j \in \ind \setminus \{0\} } I_j$ 
and, for $k \ge 2$, $T( I_{-k} ) = I_{-k+1}$.
Let $\pa_o^n := \bigvee_{k=0}^{n-1} T^{-k} \pa_o$ denote the
refinement of $\pa_o$ induced by the dynamics up to time $n$. Its 
elements are given by
\begin{equation} \label{ex-5} 
  I_{\bj^n} := I_{j_0} \cap T^{-1} I_{j_1} \cap \cdots \cap T^{-n+1} 
  I_{j_{n-1}},
\end{equation}
where $\bj^n := (j_0, \ldots, j_{n-1}) \in (\ind_o)^n$. Since $T$ is 
uniformly expanding in any given compact subset of $\R^+$ and 
clearly no \o\ converges to $+\infty$, the definition (\ref{ex-5}) implies 
that, for any infinite sequence $(j_n)_{n \in \N} \subset \ind_o$,
\begin{equation} \label{ex-10}
  \lim_{n \to \infty} \leb ( I_{(j_0, \ldots, j_{n-1})} ) = 0.
\end{equation}
Therefore, any $x$ whose forward \o\ never intersects 
$\{ a_j \}_{j \in \ind}$ (the ``boundary'' of $\pa$) has a unique 
\emph{itinerary} $(j_n)$ w.r.t.\  $\pa_o$. This means that $T^n (x) \in 
I_{j_n}$, $\forall n \in \N$; equivalently, $x \in I_{(j_0, \ldots, j_{n-1})}$, 
$\forall n \in \N$. Thus, a.e.\ $x$ has this property.

The distortion lemma that we need goes as follows:

\begin{lemma} \label{lem-dist}
  There exists $D > 1$ such that, for any $n \in \N$; any $\bj^{n+1} =
  (j_0, \ldots, j_n) \in (\ind_o)^{n+1}$ with $\leb( I_{\bj^{n+1}} ) >
  0$ and such that at least one of its components $j_k > 0$; and any
  $B \subseteq I_{\bj^{n+1}}$, one has:
  \begin{itemize}
  \item[(i)] $T^n B \subseteq I_{j_n}$;
  \item[(ii)] $\leb( T^n B \,|\, I_{j_n} ) \le D \, \leb( B \,|\,
    I_{\bj^{n+1}} )$.
  \end{itemize}
\end{lemma}

\begin{remark}
  The hypothesis that $j_k > 0$, for some $0 \le k \le n$, means that 
  the partial itinerary $\bj^{n+1}$ includes an interval $I_j$ with $j>0$. 
  This is not an unduly restrictive condition, because 
  $J = \bigcup_{j \ge 1} I_j$ is a global cross-section, so the itinerary of 
  a.e.\ $x$ will verify the hypothesis, for a large enough $n$.
\end{remark}

The statement of Lemma \ref{lem-dist} is the same as Lemma 2.3 in 
\cite{lsimple}, except that the latter has a third assertion which we 
do not need here, because we verified condition (\ref{ex-2})
by other means. The proof of Lemma \ref{lem-dist} is also
practically identical to the proof of \cite[Lem.~2.3]{lsimple}, save for
two minor changes:

\begin{enumerate}

\item For the proof of \emph{(ii)} it suffices to show 
that
\begin{equation} 
  \left| \sum_{k=0}^{n-1} \log \frac{ |T'(x_k)| } { |T'(y_k)| } \right| 
  \le C,
\end{equation}
which differs from \cite[eq.~(3.4)]{lsimple} in that $n-1$ replaces
$n$. Therefore, when parsing the \o s $(x_k)_{k=0}^{n-1}$ 
and $(y_k)_{k=0}^{n-1}$, one can posit $k_{\ell+1} := n$. Then
both $x_{k_{\ell+1}} = x_n$ and $y_{k_{\ell+1}} = y_n$ belong
in $I_{j_n}$, whence $| x_{k_{\ell+1}} - y_{k_{\ell+1}} | \le c$, 
with $c := \max_{j \in \ind _o} \leb(I_j)$. This is needed in 
\cite[eq.~(3.8)]{lsimple}. Observe that $c$ exists because, for
$j>0$, $\leb(I_j) \le a_1$ and, for $j<0$, $\leb(I_j) \le \leb(I_{-1})$.

\item Lemma 3.2 of \cite{lsimple} is replaced by 

\begin{lemma} \label{lem-my-lsy}
  There exists $C' > 0$ such that, for all $j \ge 1$, $0 \le p \le j$,
  and $x,y \in I_{-j}$,
  \begin{displaymath}
    \left| \log \frac{ (T^p)' (x) } { (T^p)' (y) } \right| \le C'
    \frac{ |T^p(x) - T^p(y)| } { L_{p-j} } \le C',
  \end{displaymath}
  where, for $p \le j-1$, $L_{p-j} := \leb(I_{p-j}) = b_{j-p} - b_{j-p-1}$
  and, for $p = j$, $L_0 := \leb(J) = a_1$ (observe that $T^p(x),
  T^p(y)$ belong to $I_{p-j}$ or $J$, respectively).
\end{lemma}

The meaning of this lemma is that the amount of distortion produced
during an `excursion' inside $I_0$ is bounded, no matter how
long the excursion. We give a detailed proof of it.

\end{enumerate}

\proofof{Lemma \ref{lem-my-lsy}} This proof is inspired by 
\cite[\S6, Lem.~5]{y}. Its main estimate, however, requires some 
original preparatory material. 

For $x \ge a_1$, set
\begin{equation} \label{my-10}
  w(x) := \int_{a_1}^x \frac1 {u(y)} \, dy.
\end{equation}
By virtue of (B4), the above defines a strictly increasing 
diverging \fn. Its inverse $v := w^{-1}$ is an increasing, concave, 
asymptotically flat bijection $[0, +\infty) \into [a_1, +\infty)$. One
verifies immediately that
\begin{align} 
\label{my-20}
  v' &= u \circ v; \\ 
\label{my-21}
  \ds \frac{v''} {v'} &= u' \circ v.
\end{align}
For $n \in \Z^+$, denote $E_n := [v(n-1),  v(n))$: these intervals
partition $[a_1, +\infty)$. For all $k \in \N$, let $n_k$ be the unique 
positive integer such that $b_k \in E_{n_k}$. We claim 
that the two partitions $\{ I_{-k} \}$ and $\{ E_n \}$ have similar 
densities. More precisely, there exists $C_1>1$ such that 
\begin{equation} \label{my-30}
  C_1^{-1} \le \frac{ \leb(I_{-k}) }{ \leb(E_{n_k}) } \le C_1 . 
\end{equation}
This entails that each interval of one partition intersects a bounded 
number of intervals of the other partition.

In fact, consider $k \ge 1$. The definitions of $b_k$ and $u$ give
\begin{equation} \label{my-40}
  \leb(I_{-k}) = b_k - b_{k-1} = b_k - T(b_k) = u(b_k).
\end{equation}
On the other hand, by the Mean Value Theorem and
(\ref{my-20}), there exists $\xi_k \in (n_k - 1, n_k)$ such that
\begin{equation} \label{my-45}
  \leb(E_{n_k}) = v(n_k) - v(n_k-1) = v'( \xi_k ) = u(v(\xi_k)).
\end{equation}
As $u$ is decreasing, both $u(b_k)$ and $u(v(\xi_k))$ lie
in the interval $( u( v(n_k) ) , u( v(n_k - 1) ] =
( v'(n_k) ) , v'(n_k - 1) ]$. Therefore, in view of 
(\ref{my-40})-(\ref{my-45}), the claim (\ref{my-30}) will be proved
if we show that
\begin{equation} \label{my-50}
  \log v'(n_k - 1) - \log v'(n_k) \le C_2,
\end{equation}
for some $C_2>0$, independent of $k$. Using again the Mean Value 
Theorem, and (\ref{my-21}), we can rewrite the above l.h.s.\ 
as $- u'(v(\eta_k))$, for some $\eta_k \in (n_k - 1, n_k)$.
But $u'$ is bounded by the assumptions on $T$, so both
(\ref{my-50}) and (\ref{my-30}) hold true.

Now for the core arguments. Take $j \ge 1$, $0 \le p \le j$,
and $x,y \in I_{-j}$, as in the statement of the lemma. 
For $0 \le i \le p-1$, there exists $\zeta_i$ 
between $T^i(x)$ and $T^i(y)$ 
(hence $\zeta_i \in I_{i-j}$) such that
\begin{equation} \label{my-70}
  \log T'( T^i(x) ) - \log T'( T^i(y) ) = \frac{ T''(\zeta_i) } 
  {T'(\zeta_i) } \left( T^i(x) - T^i(y) \right).
\end{equation}
We will estimate each term in the above r.h.s.\ separately.
To start with, $T'(\zeta_i) \ge 1$. Also, $\zeta_i \in I_{i-j}$ implies
that $\zeta_i > b_{j-i-1} \ge v(n_{j-i-1} - 1)$. The hypothesis
on $u''$, cf.\ (B4), then gives $| T''(\zeta_i) | \le u''( v(n_{j-i-1} - 1) )$.
Finally, using (\ref{my-30}), (\ref{my-45}), and the monotonicity of
$u \circ v = v'$, we obtain $| T^i(x) - T^i(y) | \le \leb(I_{i-j}) \le C_1 
u(v(\xi_{j-i})) \le C_1 v'(n_{j-i-1} - 1)$. 

All this implies that, for all $0 \le q \le p$,
\begin{equation} \label{my-80}
\begin{split}
  \left| \log \frac{ (T^q)' (x) } { (T^q)' (y) } \right| &\le 
  \sum_{i=0}^{q-1} \frac{ | T''(\zeta_i) | } {T'(\zeta_i) } \left| T^i(x) 
  - T^i(y) \right| \\
  &\le C_1 \sum_{i=0}^{q-1} u''( v(n_{j-i-1} - 1) ) \, v'(n_{j-i-1} - 1)
\end{split}
\end{equation}
Now, $(n_k)$ is an increasing, but not necessarily strictly increasing,
sequence. However, by (\ref{my-30}) \emph{et seq.}, it has bounded
multiplicity in the sense that $\# \rset{k \in \N} {n_k = j} \le C_1$.
Therefore, continuing from (\ref{my-80}),
\begin{equation} \label{my-90}
\begin{split}
  \left| \log \frac{ (T^q)' (x) } { (T^q)' (y) } \right| &< 
  C_1^2 \sum_{n=1}^\infty u''( v(n - 1) ) \, v'(n - 1) \\
  &\le C_1^2 \int_0^\infty u''(v(x)) \, v'(x) \, dx \\
  &= C_1^2 \, |u'(0)| =: C_3, 
\end{split}
\end{equation}
having used the monotonicity of $u'' \circ v$ and $v'$.

The above holds for a generic pair $x,y \in I_{-j}$, not necessarily 
the one given in the statement of the lemma. Standard arguments 
imply that
\begin{equation} \label{my-100}
  e^{-C_3} \, \frac{ |x -y| } {L_{-j}} \le \frac{ | T^q(x) - T^q(y) | } 
  {L_{q-j}} \le e^{C_3} \, \frac{ |x -y| } {L_{-j}}.
\end{equation}
Comparing the above expression for a generic $q = i \in \{ 0, \ldots,
p-1 \}$ with the same for $q = p$, we see that, for all $0 \le i \le
p-1$,
\begin{equation} \label{my-110}
  \frac{ | T^i(x) - T^i(y) | } {L_{i-j}} \le e^{2C_3} \frac{ | T^p(x) - 
  T^p(y) | } {L_{p-j}}.
\end{equation}

Using (\ref{my-110}) in the first line of (\ref{my-80}), evaluated for
$q = p$, yields
\begin{equation} \label{my-120}
\begin{split}
  \left| \log \frac{ (T^p)' (x) } { (T^p)' (y) } \right| &\le 
  e^{2C_3} \frac{ | T^p(x) - T^p(y) | } {L_{p-j}}
  \sum_{i=0}^{p-1} \frac{ | T''(\zeta_i) | } {T'(\zeta_i) } L_{i-j} \\
  &\le C' \, \frac{ | T^p(x) - T^p(y) | } { L_{p-j} }, 
\end{split}
\end{equation}
where $C' := C_3 \, e^{2C_3}$. This is so because the sum in
the first line of (\ref{my-120}) is estimated exactly in the same
way as (\ref{my-80})-(\ref{my-90}).
\qed

\section{Appendix: Proofs of technical results}
\label{sec-tech}

In this section we give the proofs of a couple of purely technical
results.

\subsection{Proof of Proposition \ref{prop-hitting}}
\label{subs-pf-hitting}

Any function $F: \ui \into \C$ defined as in (\ref{step}) is in 
$L^\infty((0,1),\mu)$. To show that it is a global \ob\ it remains to show
the existence and the value of its infinite-volume average 
\begin{equation} \label{limite}
  \avg(F) := \lim_{a\to 0^+} \frac1 {\mu([a,1))} \int_a^1 F\, d\mu.
\end{equation}

Recalling the definition of the partition $\{ B_k \}_{k \in \N}$, denote by
$( \beta_k )_{k \in \N}$ the decreasing sequence in $(0,1]$ such that 
$B_k = (\beta_{k+1}, \beta_k) \cap \rat$. Thus $\beta_0=1$.

We first take the limit (\ref{limite}) along the sequence $a =
\beta_{k+1}$. Keeping in mind that $F$ is constant on the elements of 
$\{ B_k \}$, we have
\begin{equation} \label{cla10}
  \int_{\beta_{k+1}}^1 \!\! F\, d\mu = \sum_{p=0}^k \mu(B_p) \, F|_{B_p} = 
  \sum_{j=0}^{q-1} \sum_{\pippo}^k \!\!\! f_j \: \mu(B_p)
\end{equation}
and
\begin{equation} \label{cla20}
  \mu([\beta_{k+1},1)) = \sum_{p=0}^k\, \mu(B_p).
\end{equation}
Let us introduce the notation $r_k := \mu(B_k)$ and $B_{1,k} := T^{-1} 
B_k \cap J$. Since $T$ has full branches and $\mu$ is invariant, we
see that $T^{-1} B_k = B_{k+1} \cup B_{1,k}$, with $B_{k+1} \cap 
B_{1,k} = \emptyset$, whence $\mu(B_k) = \mu(B_{k+1}) + \mu(B_{1,k})$. 
But $\mu(B_{1,k})>0$ and $\mu$ infinite, therefore the sequence $(r_k)$ is 
decreasing and the series $\sum_k r_k$ is diverging. In light of 
(\ref{limite})-(\ref{cla20}), and using the notation (\ref{mu-v-ui}), we write
\begin{equation} \label{successione}
  \lim_{k\to \infty} \mu_{[\beta_{k+1},1)} (F) = \lim_{k\to \infty} \, 
  \left( \sum_{p=0}^k r_p \right)^{\!\! -1} \: \sum_{j=0}^{q-1} \, f_j \!\!\!
  \sum_{\pippo}^k \!\!\!\ r_p 
\end{equation}
We claim that the above limit exists and equals $q^{-1} 
\sum_{j=0}^{q-1} f_j$.

Fix $j \in \{ 0,\dots, q-1 \}$ and set 
\begin{equation} \label{cla30}
  S_{j,k} := \left( \sum_{p=0}^k r_p \right)^{\!\! -1} \!\!\! \sum_{\pippo}^k \!\!\! 
  r_p \: = \: \left( \sum_{p=0}^k r_p \right)^{\!\! -1} \: \sum_{\ell=0}^{l_{j,k}} 
  r_{j+\ell q}
\end{equation}
with $l_{j,k} = \lfloor (k-j) / q \rfloor$. The claim made in the previous 
paragraph will be proved once we show that
\begin{equation} \label{cla-claim}
  \lim_{k \to \infty} S_{j,k} = \frac1q.
\end{equation}
To achieve our goal, we fix $j' \in \{ 0,\dots, q-1 \}$, with $j' \ne j$, and 
compare $S_{j,k}$ with $S_{j',k}$, for $k$ large. More in detail, we 
subdivide the finite sequence $( r_p )_{p=0}^k$ in blocks of size $q$, 
summing only the $j^\mathrm{th}$ element, respectively the 
$(j')^\mathrm{th}$ element, from each block. Upon renormalization by the 
term $\sum_{p=0}^k r_p$, we verify that the two sums have the same 
asymptotics. 

Let us implement the plan: Without loss of generality assume that $j < j'$. 
Since $( r_k )$ is decreasing, 
\begin{align}
  \sum_{\ell=0}^{l_{j,k}} r_{j+\ell q} &\ge \sum_{\ell=0}^{l_{j,k}} 
  r_{j'+\ell q} \, ; \label{cla40} \\
  \sum_{\ell=1}^{l_{j,k}} r_{j+\ell q} &\le \sum_{\ell=1}^{l_{j,k}} 
  r_{j'+ (\ell-1) q} = \sum_{\ell=0}^{l_{j,k}-1} r_{j'+\ell q} \, .
  \label{cla50} 
\end{align}
Evidently, both the above r.h.sides differ by $\sum_{\ell=0}^{l_{j',k}} 
r_{j'+\ell q}$ by a bounded quantity. Also, the l.h.s.\ of (\ref{cla50})
equals $\sum_{\ell=0}^{l_{j,k}} r_{j+\ell q} - r_j$. Dividing all terms by 
$\sum_{p=0}^k r_p$, which diverges as $k \to \infty$, we conclude
that
\begin{align}
  \limsup_{k\to \infty} S_{j,k} &= \limsup_{k\to \infty} S_{k,j'} ; \\
  \liminf_{k\to \infty} S_{j,k} &= \liminf_{k\to \infty} S_{k,j'} .
\end{align}
On the other hand, by the definition (\ref{cla30}),
\begin{equation}
  \limsup_{k\to \infty} \sum_{j=0}^{q-1}\, S_{j,k} = 
  \liminf_{k\to \infty} \sum_{j=0}^{q-1}\, S_{j,k} = 1,
\end{equation}
which implies (\ref{cla-claim}) and thus our claim.
This shows that the limit (\ref{successione}) exists and amounts 
to $q^{-1} \sum_{j=0}^{q-1} f_j$.

It remains to prove that the full limit (\ref{limite}) is the same.
For $a\in (\beta_{k+1},\beta_k)$, write
\begin{equation}
  \frac1 {\mu([a,1))} \int_a^1\, F\, d\mu = \frac1 { \mu([a,\beta_k)) + 
  \mu([\beta_k,1)) } \left( \int_{a}^{\beta_k} \! F\, d\mu + \int_{\beta_k}^1 
  F\, d\mu \right)
\end{equation}
and notice that $\mu([a,\beta_k)) \le \mu(B_k) \le \mu(B_0)$, and 
$| \int_{a}^{\beta_k} F\, d\mu | \le \|F\|_{L^\infty}\, \mu(B_0)$. Since 
$\mu([\beta_k,1))$ diverges, as $k \to \infty$, namely as $a \to 0^+$, we 
conclude that
\begin{equation}
  \lim_{a\to 0^+} \mu_{[a,1)} (F) = \lim_{k\to \infty} \mu_{[\beta_{k+1},1)} 
  (F) = \frac1q \sum_{j=0}^{q-1} f_j.
\end{equation}
The proposition is proved.
\qed

\subsection{Proof of Lemma \ref{increasing-ui}}
\label{subs-increasing-ui}

Let us fix $x,y \in (0,1)$ with $x<y$. By the Mean 
Value Theorem there exist $\xi\in (\ibr_{0}(x), \ibr_{0}(y))$ and 
$\eta \in (\ibr_{1}(y), \ibr_{1}(x))$ such that
\begin{align}
  & g(\ibr_{0}(x)) = g(\ibr_{0}(y)) + g'(\xi)\, (\ibr_{0}(x) - \ibr_{0}(y)) ;
  \label{mvt1} \\
  & g(\ibr_{1}(x)) = g(\ibr_{1}(y)) + g'(\eta)\, (\ibr_{1}(x) - \ibr_{1}(y)) .
  \label{mvt2}
\end{align}
Using (\ref{f1}), (H1), the identity $\bb = 1-\aa$ 
and the inequality $g'(\xi) \ge g'(\eta)$, which follows from (H2) and 
the concavity of $g$, we can write:
\begin{align}
   g_{1}(y) - g_{1}(x) &= \aa(y) g(\ibr_{0}(y)) + \bb(y) g(\ibr_{1}(y)) - 
   \aa(x) g(\ibr_{0}(y)) - \bb(x) g(\ibr_{1}(y)) \nonumber \\[7pt]
   &\qquad - \aa(x) g'(\xi) (\ibr_{0}(x) - \ibr_{0}(y)) - \bb(x) g'(\eta) 
   (\ibr_{1}(x) - \ibr_{1}(y)) \nonumber \\[4pt]
   &\ge \aa(y) \Big( g(\ibr_{0}(y)) - g(\ibr_{1}(y)) \Big) - \aa(x) \Big( 
   g(\ibr_{0}(y)) - g(\ibr_{1}(y)) \Big) \nonumber \\
   &\qquad - g'(\eta) \Big( (\aa(x)-\bb(x)) (\ibr_{0}(x)- \ibr_{0}(y)) \\
   &\qquad + \bb(x) (\ibr_{0}(x)+ \ibr_{1}(x)- \ibr_{0}(y)- \ibr_{1}(y))\Big) .
   \nonumber
\end{align}
having also used that $\ibr_0$ is increasing. We study the last 
term in the above inequality piece by piece. By (H2) and the 
monotonicity of $g$; the first assertion of (H7); (H6); (H1); (H3), we obtain,
respectively:
\begin{align}
  g(\ibr_{0}(y)) - g(\ibr_{1}(y)) &\le 0 ; \\
  \aa(y) - \aa(x) &\le 0 ; \\
  \aa(x)-\bb(x) &\ge 0 ; \\
  \ibr_{0}(x)- \ibr_{0}(y) &\ge 0 ; \\
  \ibr_{0}(x)+ \ibr_{1}(x)- \ibr_{0}(y)- \ibr_{1}(y) &\le 0 .
\end{align}
Hence $g_{1}(y)-g_{1}(x) \ge 0$, proving that $g_1$ is increasing.

\skippar

We now show that $g_1$ is concave, namely, for any pair $x,y \in \ui$, 
$x<y$, and $z = tx + (1-t) y \in (x,y)$, with $0 \le t \le 1$, we verify that
\begin{equation} \label{cla80}
  g_{1}(z) \ge t\, g_{1}(x) + (1-t)\, g_{1}(y).
\end{equation}

Clearly $g_{1}(z) = t\, g_{1}(z) + (1-t)\, g_{1}(z)$. By means of (\ref{f1})
we have
\begin{align}
  g_{1}(z) &- t\, g_{1}(x) - (1-t)\, g_{1}(y) \nonumber \\[4pt]
  &=  t \Big( \aa(z) g(\ibr_{0}(z)) + \bb(z) g(\ibr_{1}(z)) - \aa(x) 
  g(\ibr_{0}(x)) - \bb(x) g(\ibr_{1}(x)) \Big) \\
  &\qquad + (1-t) \Big( \aa(z) g(\ibr_{0}(z)) + \bb(z) g(\ibr_{1}(z)) - 
  \aa(y) g(\ibr_{0}(y)) - \bb(y) g(\ibr_{1}(y)) \Big). \nonumber
\end{align}
We apply the Mean Value Theorem, as in (\ref{mvt1})- (\ref{mvt2}) to
write:
\begin{align}
  & g(\ibr_{0}(x)) = g(\ibr_{0}(z)) + g'(\xi_{1})\, (\ibr_{0}(x) - \ibr_{0}(z)) ;
  \label{mvt3} \\
  & g(\ibr_{1}(x)) = g(\ibr_{1}(z)) + g'(\eta_{1})\, (\ibr_{1}(x) - \ibr_{1}(z)) ;
  \label{mvt4} \\
  & g(\ibr_{0}(y)) = g(\ibr_{0}(z)) + g'(\xi_{2})\, (\ibr_{0}(y) - \ibr_{0}(z)) ;
  \label{mvt5} \\
  & g(\ibr_{1}(y)) = g(\ibr_{1}(z)) + g'(\eta_{2})\, (\ibr_{1}(y) - \ibr_{1}(z)) , 
  \label{mvt6}
\end{align}
for some $\xi_{1}\in (\ibr_{0}(x), \ibr_{0}(z))$,  $\xi_{2}\in 
(\ibr_{0}(z), \ibr_{0}(y))$, $\eta_{1}\in (\ibr_{1}(z), \ibr_{1}(x))$ and 
$\eta_{2}\in (\ibr_{1}(y), \ibr_{1}(z))$.
Making the substitutions (\ref{mvt3})-{\ref{mvt6}) yields
\begin{align}
  g_{1}(z) &- t\, g_{1}(x) - (1-t)\, g_{1}(y) \nonumber \\[4pt]
  &= t \Big[ \Big( \aa(z) g(\ibr_{0}(z)) + \bb(z) g(\ibr_{1}(z)) \Big) \! - \! 
  \Big( \aa(x) g(\ibr_{0}(z)) + \bb(x) g(\ibr_{1}(z)) \Big) \Big] \nonumber \\
  &\qquad +(1-t) \Big[ \Big( \aa(z) g(\ibr_{0}(z)) + \bb(z) g(\ibr_{1}(z)) 
  \Big) \! - \! \Big( \aa(y) g(\ibr_{0}(z)) + \bb(y) g(\ibr_{1}(z)) \Big) \Big] 
  \nonumber \\
  &\qquad -\Big[ t \Big( \aa(x) g'(\xi_{1}) \, (\ibr_{0}(x) - \ibr_{0}(z)) + 
  \bb(x) g'(\eta_{1})\, (\ibr_{1}(x) - \ibr_{1}(z)) \Big) \label{cla90} \\
  &\qquad + (1-t) \Big( \aa(y) g'(\xi_{2})\, (\ibr_{0}(y) - \ibr_{0}(z)) + 
  \bb(y) g'(\eta_{2})\, (\ibr_{1}(y) - \ibr_{1}(z)) \Big) \Big] \nonumber 
  \\[1pt]
  &=: \Theta_1 - \Theta_2, \nonumber
\end{align}
where $\Theta_1$ corresponds the second and third lines above, and 
$\Theta_2$ to the opposite of the fourth and fifth lines, cf.\ (\ref{cla100}) 
and (\ref{cla110}) below.

Let us first consider $\Theta_1$. Since $\bb=1-\aa$, we can write
\begin{align}
  \Theta_1 &:= t \Big[ \Big( \aa(z) g(\ibr_{0}(z)) + \bb(z) g(\ibr_{1}(z)) \Big) - 
  \Big( \aa(x) g(\ibr_{0}(z)) + \bb(x) g(\ibr_{1}(z)) \Big) \Big] \nonumber \\
  &\qquad +(1-t) \Big[ \Big( \aa(z) g(\ibr_{0}(z)) + \bb(z) g(\ibr_{1}(z)) \Big) - 
  \Big( \aa(y) g(\ibr_{0}(z)) + \bb(y) g(\ibr_{1}(z)) \Big) \Big] \nonumber 
  \\[4pt]
  &= \aa(z) g(\ibr_{0}(z)) + \bb(z) g(\ibr_{1}(z)) \label{cla100} \\[4pt]
  &\qquad -\Big[ \Big( t\aa(x)+(1-t)\aa(y)\Big) g(\ibr_{0}(z)) +  \Big( 
  t\bb(x)+(1-t)\bb(y)\Big) g(\ibr_{1}(z))\Big] \nonumber \\[2pt]
  & = \aa(z) \Big( g(\ibr_{0}(z)) - g(\ibr_{1}(z))\Big) - \Big( 
  t\aa(x)+(1-t)\aa(y)\Big) \Big( g(\ibr_{0}(z)) - g(\ibr_{1}(z))\Big) 
  \nonumber \\
  & = \Big( \aa(z) - t\aa(x)-(1-t)\aa(y) \Big) \, \Big( g(\ibr_{0}(z)) - 
  g(\ibr_{1}(z))\Big). \nonumber
\end{align}
On the other hand, using the convexity of $\aa$, cf.\ (H7), the 
monotonicity of $g$ and (H2) we obtain:
\begin{align}
  \aa(z) - t\aa(x)-(1-t)\aa(y) &\le 0; \\
  g(\ibr_{0}(z)) - g(\ibr_{1}(z)) &\le 0.
\end{align}
Hence $\Theta_1 \ge 0$. 

As for $\Theta_2$, we recall the definitions of $\xi_i, \eta_i$, 
$i \in \{ 1,2 \}$, given in (\ref{mvt3})-(\ref{mvt6}). Since $\ibr_0$ is 
increasing and $\ibr_1$ is decreasing, we have the ordering 
$\eta_1 > \eta_2 > \xi_2 > \xi_1$, whence
\begin{equation} \label{f1-conc}
  g'(\xi_1) \ge g'(\xi_2) \ge g'(\eta_2) \ge g'(\eta_1),
\end{equation}
because $g$ is concave. Now, by (H7) and (H5), $\aa$ is decreasing 
and $\bb$ is increasing. Therefore:
\begin{equation} \label{cla110}
\begin{split}
  \Theta_2 &:= t \Big( \aa(x) g'(\xi_{1}) \, (\ibr_{0}(x) - \ibr_{0}(z)) + 
  \bb(x) g'(\eta_{1})\, (\ibr_{1}(x) - \ibr_{1}(z)) \Big) \\
  &\qquad +(1-t) \Big( \aa(y) g'(\xi_{2})\, (\ibr_{0}(y) - \ibr_{0}(z)) + 
  \bb(y) g'(\eta_{2})\, (\ibr_{1}(y) - \ibr_{1}(z)) \Big) \\
  &\le t \Big( \aa(y) g'(\xi_{2}) \, (\ibr_{0}(x) - \ibr_{0}(z)) + \bb(y) 
  g'(\eta_{2})\, (\ibr_{1}(x) - \ibr_{1}(z)) \Big) \\
  &\qquad +(1-t) \Big( \aa(y) g'(\xi_{2})\, (\ibr_{0}(y) - \ibr_{0}(z)) + 
  \bb(y) g'(\eta_{2})\, (\ibr_{1}(y) - \ibr_{1}(z)) \Big) \\
  &= g'(\xi_{2}) \aa(y) \Big( t \ibr_{0}(x) + (1-t) \ibr_{0}(y) - \ibr_{0}(z)
  \Big) \\
  &\qquad +g'(\eta_{2}) \bb(y) \Big( t \ibr_{1}(x) + (1-t) \ibr_{1}(y) - 
  \ibr_{1}(z)\Big).
\end{split}
\end{equation}
Using the concavity of $\ibr_0$, cf.\ (H4), and (\ref{f1-conc}), we have:
\begin{align}
  \Theta_2 &\le  g'(\eta_{2}) \Big[ \aa(y) \Big( t \ibr_{0}(x) + 
  (1-t) \ibr_{0}(y) - \ibr_{0}(z)\Big) \nonumber \\
  &\qquad +\bb(y) \Big( t \ibr_{1}(x) + (1-t) \ibr_{1}(y) - \ibr_{1}(z)\Big) 
  \Big] \label{cla120} \\
  &= g'(\eta_{2}) \Big[ (\aa(y)-\bb(y)) \Big( t \ibr_{0}(x) + (1-t) 
  \ibr_{0}(y) - \ibr_{0}(z)\Big) \nonumber \\
  &\qquad +\bb(y) \Big( t (\ibr_{0}(x)+ \ibr_{1}(x)) + (1-t) (\ibr_{0}(y) + 
  \ibr_{1}(y)) - (\ibr_{0}(z)+ \ibr_{1}(z))\Big) \Big] . \nonumber
\end{align}
Now, by the hypotheses on $g$ and $\bb$, $g'(\eta_{2}) \ge 0$, 
$\bb(y) > 0$. By (H6), $\aa(y)-\bb(y) \ge 0$. Moreover, (H4) gives:
\begin{align}
  t \ibr_{0}(x) + (1-t) \ibr_{0}(y) - \ibr_{0}(z) &\le 0; \\
  t (\ibr_{0}(x)+ \ibr_{1}(x)) + (1-t) (\ibr_{0}(y)+ \ibr_{1}(y)) - (\ibr_{0}(z)+ 
  \ibr_{1}(z)) &\le 0.
\end{align}
Applying all these inequalities to (\ref{cla120}) shows that $\Theta_2 
\le 0$. 

Therefore $\Theta_1 - \Theta_2 \ge 0$, which, in view of (\ref{cla90}), 
proves our claim (\ref{cla80}).
\qed

\end{document}